\documentclass[aap,MSNbibl,citesort,seceqn,dvips]{arximspdf}

\doi{10.1214/12-AAP844} 
\volume{23}
\issue{1}
\pubyear{2013}
\firstpage{308}
\lastpage{347}

\makeatletter
\newtheorem{thmm}{Theorem}[section]
\newtheorem{lem}[thmm]{Lemma}

\newtheorem{prop}[thmm]{Proposition}
\newproclaim{rem}[thmm]{Remark}
\newproclaim{eg}[thmm]{Example}
\newproclaim{defn}[thmm]{Definition}
\newproclaim{assum}[thmm]{Assumption}
\newproclaim{definition}{Definition}

\newproclaim{remark}{Remark}
\newproclaim{example}{Example}

\newproclaim{assumption}{Assumption}

\newcommand{\hP}{{\hat\dbP}}
\newcommand{\eps}{\varepsilon}

\def\dbD{\mathbb{D}}
\def\dbE{\mathbb{E}}
\def\dbF{\mathbb{F}}
\def\dbG{\mathbb{G}}
\def\dbH{\mathbb{H}}
\def\dbI{\mathbb{I}}
\def\dbL{\mathbb{L}}
\def\dbP{\mathbb{P}}
\def\dbR{\mathbb{R}}
\def\dbS{\mathbb{S}}
\def\dbQ{\mathbb{Q}}

\def\Dom{\operatorname{dom}}

\def\a{\alpha}
\def\b{\beta}
\def\d{\delta}
\def\e{\varepsilon}
\def\k{\kappa}
\def\n{\nu}
\def\si{\sigma}
\def\t{\tau}
\def\f{\varphi}
\def\G{\Gamma}
\def\L{\Lambda}
\def\O{\Omega}
\def\cA{\mathcal{A}}
\def\cF{\mathcal{F}}
\def\cM{\mathcal{M}}
\def\cP{\mathcal{P}}
\def\cS{\mathcal{S}}
\def\cV{\mathcal{V}}
\def\cY{\mathcal{Y}}
\def\cZ{\mathcal{Z}}
\def\cd{\cdot}
\def\cds{\cdots}

\newcommand{\dfnn}{\stackrel{\triangle}{=}}

\def\limsup{\overline{\operatorname{lim}}}
\def\esup{\operatorname{ess\,sup}}
\def\cd{\cdot}
\def\cds{\cdots}

\def\1{\mathbf{1}}
\def\SM{{\cS\cM}^2}
\def\hSM{\widehat{\cS\cM}^2}

\def\:{ \dvtx}
\makeatother

\begin{document}
\begin{frontmatter}

\title{Dual formulation of second order target problems}
\runtitle{Dual of second order target problems}

\begin{aug}
\author[A]{\fnms{H.~Mete} \snm{Soner}\ead[label=e1]{hmsoner@ethz.ch}\thanksref{t1}},
\author[B]{\fnms{Nizar} \snm{Touzi}\ead[label=e2]{nizar.touzi@polytechnique.edu}\thanksref{t2}}
\and
\author[C]{\fnms{Jianfeng} \snm{Zhang}\corref{}\ead[label=e3]{jianfenz@usc.edu}\thanksref{t3}}
\thankstext{t1}{Supported in part by the
European Research Council under the Grant 228053-FiRM and
supported by the ETH Foundation.}
\thankstext{t2}{Supported by the Chair \textit{Financial Risks} of
the \textit{Risk Foundation} sponsored by Soci\'et\'e
G\'en\'erale, the Chair \textit{Derivatives of the Future} sponsored
by the {F\'ed\'eration Bancaire Fran\c{c}aise}, and
the Chair \textit{Finance and Sustainable Development} sponsored by
EDF and Calyon.}
\thankstext{t3}{Supported in part by NSF Grants DMS-06-31366 and DMS-10-08873.}
\runauthor{H. M. Soner, N. Touzi and J. Zhang}
\affiliation{ETH Z\"urich and Swiss Finance Institute, CMAP, Ecole Polytechnique, and~University of
Southern California}
\address[A]{H. M. Soner\\
ETH (Swiss Federal Institute of Technology)\\
Z\"urich\\
and\\
Swiss Finance Institute\\
HG G 54.3, Ramistrasse 101\\
8092 Zurich\\
Switzerland\\
\printead{e1}} 
\address[B]{N. Touzi\\
CMAP, Ecole Polytechnique\\
Route De Saclay\\
91128 Palaiseau\\
France\\
\printead{e2}}
\address[C]{J. Zhang\\
Department of Mathematics\\
University of Southern California\\
3620 S. Vermont Ave, KAP 108\\
Los Angeles, California 90089\\
USA\\
\printead{e3}}
\end{aug}

\received{\smonth{4} \syear{2010}}
\revised{\smonth{1} \syear{2012}}

\begin{abstract}
This paper provides a new formulation of second order
stochastic target problems introduced in [\textit{SIAM J. Control Optim.} \textbf{48} (2009) 2344--2365]
by modifying the
reference probability so as to allow for different scales.
This new ingredient enables us to
prove a dual formulation of the target problem as the supremum
of the solutions of standard backward stochastic differential
equations. In particular, in the Markov case, the dual problem is
known to be connected to a fully nonlinear, parabolic
partial differential equation and this
connection can be viewed as a stochastic representation
for all nonlinear, scalar, second order, parabolic equations
with a convex Hessian dependence.
\end{abstract}

%

\end{frontmatter}

\section{Introduction}


The connection between the backward stochastic differential equations
(BSDE hereafter) and the nonlinear, parabolic partial
differential equations (PDE hereafter) is well documented.
Indeed, the standard BSDEs,
as introduced by Pardoux and Peng~\cite{pardouxpeng},
are known to provide a stochastic
representation for the solutions of semi-linear PDEs
in the Markov case.
In this representation,
the diffusion coefficient of the underlying process
is the linear coefficient of the Hessian variable in the PDE.
Therefore, the connection
to fully nonlinear equations requires an extension
that should allow for stochastic processes
with different diffusion coefficients.
Indeed,~\cite{cstv} develops
such a generalization to the second order
and also proves a Markovian
uniqueness result in an appropriate class.
However, no existence theory is available for this
generalization, with the one exception in the Markov context.
In this case any smooth solution of the related PDE, if it exists,
is easily seen to be a solution of the second order BSDE.
A closely related class of control problems,
called the \textit{second order stochastic target problem}, was
introduced in~\cite{ST} as well.

In this paper we provide a new formulation
for the second order stochastic target problems.
A better understanding of the target problem is
essential for a coherent theory of second order BSDEs.
Indeed, we develop this theory in
our accompanying work~\cite{STZ09d},
including existence and uniqueness results
with minimal assumptions.

We continue with the description of
the target problem. Let $B$ be a
Brownian motion under the probability measure
$\dbP_0$ and $\{\cF_t,t\ge0\}$
be the corresponding filtration.
For a continuous semimartingale $Z$,
we denote by $\Gamma$ the density of its
covariation with $B$. We then define the
controlled process $Y$ by
%
\begin{equation}\label{YyZ}
\qquad Y_t := y - \int_0^t H_s(Y_s,Z_s,\Gamma_s)\,ds
+ \int_0^t Z_s\circ dB_s,
\qquad
d\langle Z,B\rangle_t = \Gamma_t \,dt,
\end{equation}
where $\circ$ denotes the Fisk--Stratonovich
stochastic integration.
We assume that the given random nonlinear function
$H$ satisfies the standard
Lipschitz and measurability conditions.
Then, for any reasonable process $Z$ and an
initial condition $y$, a unique solution,
which is denoted by $Y^{y,Z}$, exists.
We now fix a time
horizon, say, $T=1$, and a class
of admissible controls $\cZ^0$. Then, given an
$\cF_1$ measurable random variable $\xi$,
\cite{ST} defines the
second order stochastic target problem by
%
\begin{equation}\label{target0}
\cV^0 := \inf\{y\dvtx  Y^{y,Z}_1\ge\xi\
\dbP_0\mbox{-a.s. for some } Z\in\cZ^0 \}.
\end{equation}
In this formulation, the structure of the set of
admissible controls
is crucial. In fact, if $\cZ^0$ is not properly
defined, then the dependence of the problem
on the variable $\Gamma$
can be trivialized. We refer to~\cite{cetin-s-t}
for a detailed discussion of this issue in
a particular example of mathematical finance.
One of the achievements of the approach
given below is to avoid this strong dependence
on the control set and
simply to work with standard spaces.

As in many optimization problems,
convex duality results provide a deeper
understanding and powerful technical tools.
Indeed, they are an essential
step for the well-posedness of the
second order backward stochastic differential equations,
as proved in our accompanying paper~\cite{STZ09d}.
Motivated by these, we adopt a new point of view
for the target problems
which also allows for the construction of the dual.
This new formulation differs
from that of~\cite{ST} in two instances.
First, we reinforce the constraint
$Y^{y,Z}_1\ge\xi$ in (\ref{target0}) by requiring
that it should hold under various mutually singular measures
and not only on the support of $\dbP_0$. Second,
the set of admissible controls utilized here
is more natural and, as discussed above,
it avoids the technical aspects of~\cite{ST}.

Our reformulation is motivated by the work
of Denis and Martini~\cite{DM}
on the deep theory of quasi-sure stochastic analysis.
An important related probabilistic
notion, introduced by Peng~\cite{Peng-G},
is the $G$-Brownian motion.
Here instead of using these two powerful
tools, we employ a direct approach by assuming
sufficient regularity.
One drawback of all these approaches is the
implicit regularity assumption.
Indeed, in all these approaches, integrability
in any power is possible only if the random variable
is quasi-surely continuous. This is a Lusin type
of result and is not restrictive when there is only
countably many measures. However, in general,
this is an additional constraint. In one of our accompanying
papers~\cite{STZ09a}, we provide an alternative approach
through aggregation of random variables.
The general aggregation result of~\cite{STZ09a}
allows us to consider a larger class of random variables,
but then the class of probability measures must be slightly restricted.

We believe our approach has several advantages:
\begin{itemize}[--]
\item[--] It avoids redeveloping an appropriate theory of stochastic integration
from scratch, as it is done in~\cite{DM} and~\cite{Peng-G}.

\item[--] More importantly, a representation theorem
is available in our framework as proved in~\cite{STZ09b}.

\item[--] Finally, by deriving appropriate estimates,
it is shown in~\cite{STZ09d} that one can
extend these concepts to a larger space
with regularity conditions.
Indeed, a similar extension of $G$-martingales
is given in~\cite{denishupeng},
showing that they
cover the same space as in the quasi-sure analysis of~\cite{DM}.
\end{itemize}

We next provide an intuitive description
of our formulation. For this heuristic
explanation we assume a Markov structure.
Namely, we assume that $H$ in (\ref{YyZ})
and $\xi$ in (\ref{target0}) are given by
%
\begin{equation}
\label{Markov}
H_t(y,z,\gamma) = h(t,X_t,y,z,\gamma),
\qquad
\xi=g(X_T),
\end{equation}
where $dX_t = dB_t$
and $h$, $g$ are deterministic scalar functions.
Let $\cV^0(t,x)$ be defined as in (\ref{target0})
with time origin at $t$ and $X_t=x$.
As it is usual, we assume that
$\gamma\mapsto h(t,x,y,z,\gamma)$ is nondecreasing.
Then, by an appropriate choice of admissible controls $\cZ$,
it is shown in~\cite{ST} that this problem is a viscosity
solution of the corresponding dynamic programming equation,
%
\begin{equation}\label{pdeu}
\qquad-\frac{\partial u}{\partial t} - h
(t,x,u(t,x),Du(t,x),D^2u(t,x))=0,\qquad
u(1,x)=g(x).
\end{equation}
We further assume that $\gamma\mapsto h(t,x,r,p,\gamma)$
is convex. Then,
%
\begin{equation}\label{hmarkov}
h(t,x,r,p,\gamma)
=
\sup_{a\ge0} \biggl\{\frac12 a\gamma- f(t,x,r,p,a
)\biggr\} ,
\end{equation}
where $f$ is the (partial) convex conjugate of
$h$ with respect to $\gamma$.
Let $D_f$ be the domain of $f$ as a function of $a$.
By the classical maximum principle
of parabolic differential equations,
we expect that, for every $a\in D_f$,
the solution $u\ge u^a$, where~$u$ solves (\ref{pdeu})
and $u^a$ is defined as the solution of the
following semi-linear PDE:
%
\begin{eqnarray}\label{pdeua}
-\frac{\partial u}{\partial t} - \frac12 aD^2u(t,x)
+ f(t,x,u(t,x),Du(t,x),a)=0,
\nonumber
\\[-8pt]
\\[-8pt]
\eqntext{u(1,x)=g(x).}
\end{eqnarray}
In turn, by standard results, $u^a(t,x)=Y^a_t$,
where, for $s\in[t,T]$,
%
\begin{eqnarray}
\label{FBSDEa}
X^a_s& =& x+\int_t^s a_r^{1/2} \,dB_r,
\nonumber
\\[-8pt]
\\[-8pt]
\nonumber
Y^a_s
&=&
g(X^a_T)-\int_s^T f(r,X^a_r,Y^a_r,Z^a_r,a)\,dr
-\int_t^T Z^a_r a^{1/2}\,dB_s.
\end{eqnarray}
We have formally argued that $\cV^0(t,x) \ge Y^a_t$
for any $a\in D_f$.
Let $\cA^f$ be the collection of all processes with values in $D_f$.
By extending (\ref{FBSDEa}) to processes $a$, it is then natural to
consider the problem
%
\begin{equation}\label{dualmarkov}
V_t := \sup_{a\in\cA^f} Y_t^a
\end{equation}
as the dual of the primal
stochastic target problem.
Indeed, the optimization problem (\ref{dualmarkov})
corresponds to the dual formulation of the second
order target problem in the Markov case.
Such a duality relation was suggested in the specific
example of~\cite{STgamma}
and can be proved rigorously by showing that $v(t,x):=V_t$
is a viscosity solution of the fully nonlinear PDE (\ref{pdeu}).
This, by uniqueness, implies that $v=\cV^0$. Of course, such an
argument requires some technical conditions
at least to guarantee that comparison of viscosity
supersolutions and subsolutions holds true for
the PDE (\ref{pdeu}).

The main object of this paper is to provide a
purely probabilistic proof of this duality result.
Moreover, our duality result does not require
to restrict the problem to the Markov framework.

We should mention that we use weak formulation in our approach, that
is, instead of controlling the state process~$X^a$ in (\ref{FBSDEa}),
our control is the distribution of~$X^a$ on its canonical space. See
(\ref{Xa}) below for the precise definition. Such weak formulation is
important for modeling model uncertainty, as in~\cite{DM} and \cite
{Peng-G}. In the contexts of stochastic control, which naturally uses
strong formulation, some ideas have already appeared in the literature;
see, for example, El Karoui and Quenez~\cite{EQ} and Peng \cite
{Peng2}. In particular,~\cite{Peng2} uses the notion of r.c.p.d. which
turns out to be crucial in our approach.

This paper is organized as follows. After introducing the
probabilistic structures in the next section, we provide
the definition of the stochastic target problem in
Section~\ref{sect-2target}. Two relaxations,
which are also shown to be equivalent to the original
problem, are also introduced in that section. The main
duality result is stated and proved in the following section.
Section~\ref{sectweak2target} is devoted
to a weaker formulation. An extension is outlined
in the next section and in the \hyperref[sectappendix]{Appendix} we
provide the proofs of two technical results.

\section{The setup}
\label{sect-setup}

Let $\O:= \{\omega\in C([0,1], \dbR^d)\dvtx  \omega_0=0\}$ be the
canonical space equipped with the uniform norm $\|\omega\|_\infty:=
\sup_{0\le t\le1}|\omega_t|$, $B$ the canonical process, $\dbP_0$
the Wiener measure, $\dbF:= \{\cF_t\}_{0\le t\le1}$ the filtration
generated by $B$, and $\dbF^+:=\{\cF^+_t,0\le t\le1\}$ the right
limit of $\dbF$.

We say a probability measure $\dbP$ is a local martingale measure if
the canonical process $B$ is a local martingale under $\dbP$. By F\"
{o}llmer~\cite{follmer} (see also Karandikar~\cite{Karandikar} for a
more general result), there exists an $\dbF$-progressively measurable
process, denoted as $\int_0^t B_s \,dB_s$, which coincides with It\^o's
integral, $\dbP$-a.s. for all local martingale measures $\dbP$. In
particular, this provides a pathwise definition of
\[
\langle B\rangle_t := B_t B_t^{\mathrm{T}} - 2 \int_0^t B_s \,dB^{\mathrm{T}}_s
\quad\mbox
{and}\quad \hat a_t:= \mathop{\limsup}_{\e\downarrow0} \frac{1}{\e}(\langle
B\rangle_t-\langle B\rangle_{t-\e}),
\]
where $^{\mathrm{T}}$ denotes the transposition, and the $\limsup$ is taken
componentwise and pointwise in $\omega$.
Clearly, $\langle B\rangle$ coincides with the $\dbP$-quadratic
variation of
$B$, $\dbP$-a.s. for all local martingale measures $\dbP$.

Let $\overline{\cP}_W$ denote the set of all local martingale
measures $\dbP$ such that
%
\begin{eqnarray}
\label{overlinecPW}
\langle B\rangle_t \mbox{ is absolutely continuous in } t
\mbox{ and }
\hat a \mbox{ takes values in } \dbS^{>0}_d,\quad \dbP\mbox{-a.s.},\hspace*{-35pt}
\end{eqnarray}
where $\dbS^{>0}_d$ denotes the space of all $d\times d$ real-valued
positive definite matrices. We note that, for different $\dbP_1, \dbP
_2\in\overline{\cP}_W$, in general $\dbP_1$ and $\dbP_2$ are
mutually singular. For any $\dbP\in\overline{\cP}_W$, it follows
from the L\'evy characterization that It\^o's stochastic integral under
$\dbP$,
%
\begin{equation}
\label{WP}
W^\dbP_t := \int_0^t \hat a^{-{1\slash2}}_s\,dB_s, \qquad t\in[0,1],\dbP
\mbox{-a.s.}
\end{equation}
defines a $\dbP$-Brownian motion.
As in~\cite{STZ09a}, we abuse the terminology of Denis and Martini
\cite{DM} as follows:
%
\begin{defn}
\label{defn-qs}
For any subset $\cP\subset\overline{\cP}_W$,
we say a property holds $\cP$-quasi-surely
($\cP$-q.s. for short) if it holds
$\dbP$-a.s. for all $\dbP\in\cP$.
\end{defn}

In this paper we concentrate on the subclass $\overline\cP_S \subset
\overline\cP_W$ consisting of all probability measures
%
\begin{eqnarray}
\label{Xa}
\qquad\dbP^\a:= \dbP_0 \circ(X^\a)^{-1} \qquad\mbox{where } X^\a_t := \int
_0^t \a_s^{1\slash2} \,dB_s, t\in[0,1], \dbP_0\mbox{-a.s.}
\end{eqnarray}
for some $\dbF$-progressively measurable process $\a$ taking values
in $\dbS^{>0}_d$ with
$\int_0^1 |\a_t| \,dt <\infty$, $\dbP_0$-a.s.
We recall from~\cite{STZ09a} that
%
\begin{equation}
\label{overlinecPS}
\overline{\cP}_S = \{\dbP\in\overline{\cP}_W\dvtx
\overline{\dbF^{W^\dbP}}^\dbP= \overline{\dbF}^\dbP\},
\end{equation}
where $\overline{\dbF}^\dbP$ (resp.,
$\overline{\dbF^{W^\dbP}}^\dbP$) is the $\dbP$-augmentation
of the filtration generated by $B$ (resp., by $W^\dbP$). Moreover,
%
\begin{equation}
\begin{tabular}{p{270pt}@{}}
\label{0-1 MRP}
every $\dbP\in\overline\cP_S$ satisfies
the Blumenthal zero--one law
and the martingale representation property.
\end{tabular}
\end{equation}

Notice that an $\dbF$-progressively measurable process can be viewed
as a mapping from $[0,T]\times\O$ to $\dbR^d$. Moreover, $X^\a$
takes values in $\O$ and, thus, its canonical space is also $\O$ and
the canonical filtration is still $\dbF$. We have the following simple lemma.
%
\begin{lem}
\label{lem-atilde}
Let $\a$ be an $\dbF$-progressively measurable
process taking values in $\dbS^{>0}_d$ with
$\int_0^1 |\a_t| \,dt <\infty$, $\dbP_0$-a.s. Then there exists an
$\dbF$-progressively measurable mapping
$\b_\a\dvtx  [0,T]\times\O\to\dbR^d$ such that
\begin{eqnarray*}
B& =& \b_\a(X^\a),\qquad \dbP_0\mbox{-a.s.}\quad
\mbox{and}\\
 W^{\dbP^\a} &=& \b_\a(B),\qquad
\hat a (B)= \a\circ\b_\a(B),\qquad dt\times\dbP^\a\mbox{-a.s.}
\end{eqnarray*}
\end{lem}

\begin{pf}
First, by~\cite{STZ09a}, Lemma 8.1, we know $\overline{\dbF
^{X^\a}}^{\dbP_0} = \overline{\dbF^B}^{\dbP_0}$, and, in
particular, $B$ is $\overline{\dbF^{X^\a}}^{\dbP_0}$ progressively
measurable. By~\cite{STZ09a}, Lemma 2.4 and Remark~\ref{rem-version}
below, there exists an $\dbF^{X^\a}$-progressively measurable process
$\tilde B$ such that $\tilde B = B$, $\dbP_0$-a.s. Then, by viewing
$\O$ as the canonical space of $X^\a$, one may identify the process~$\tilde B$
as an $\dbF$-progressively measurable mapping $\b_\a$.
Changing back to the canonical space of $B$ and noting that $X^\a$
takes values in $\O$, we have $\tilde B(\omega) = \b_\a(X^\a
(\omega))$ for all $\omega\in\O$, and, therefore, $B = \b_\a
(X^\a)$, $\dbP_0$-a.s.

Now it follows from the definition of $\dbP^{\hspace*{1pt}\a}$ that
%
\begin{equation}
\label{P0Pa}
(B, \tilde W^\a)_{\dbP^\a}
= (X^\a, B)_{\dbP_0} \qquad\mbox{where }
\tilde W^\a:= \b_\a(B),
\end{equation}
that is, the $\dbP^\a$-distribution of
$(B, \tilde W^\a)$ is equal to the
$\dbP_0$-distribution of $(X^\a, B)$.
Note that
$d\langle B\rangle_t = \hat a_t(B) \,dt$,
$\dbP^\a$-a.s. and $d\langle X^\a\rangle_t = \a(B)\,dt
= \a\circ\b_\a(X^\a)\,dt$, $\dbP_0$-a.s. Then
\[
(B, \tilde W^\a, \hat a(B))_{\dbP^\a}
= \bigl(X^\a, B, \a\circ\b_\a(X^\a)\bigr)_{\dbP_0}.
\]
This implies that $\hat a(B)
= \a\circ\b_\a(B)$, $dt\times\dbP^\a$-a.s.
Moreover, since $d B_t = \a^{-1\slash2}_t(B)\,dX^\a_t
= \a^{-1\slash2}_t(\b(X^\a))\,dX^\a_t$,
$\dbP_0$-a.s. it follows from (\ref{P0Pa}) that
\begin{eqnarray*}
\tilde W^\a_t = \int_0^t \a^{-1\slash2}_s(\b(B))\,dB_s
= \int_0^t \hat a^{-1\slash2}_s(B) \,dB_s = W^{\dbP^\a}_t,\qquad
t\in[0,1], \dbP^\a\mbox{-a.s.}\\[-10pt]\qed
\end{eqnarray*}
\noqed\end{pf}

\begin{rem}
\label{rem-version}
In the standard stochastic analysis literature,
the theory is developed under the augmented filtration.
Because we are working under mutually singular measures,
unless otherwise stated, we shall use the filtration $\dbF$.
We recall from~\cite{STZ09a} that, for every probability
measure $\dbP$, every $\overline\dbF^\dbP$-progressively
measurable process $X$ has an $\dbF$-progressively
measurable version $\tilde X$, that is, $X=\tilde X$, \mbox{$\dbP$-a.s.}
Therefore, given $\dbP$, all processes involved
in this paper will be considered in their
$\dbF$-version. However, notice that such a
version may depend on $\dbP$. See also Remark~\ref{rem-FF+} below.

Moreover, following similar arguments, the above result still holds
true if we replace $\dbF$ by an arbitrary filtration. In the proof of
Lemma~\ref{lem-atilde}, we have used the result on the filtration
$\dbF^{X^\a}$.
\end{rem}

Finally, we clarify that
by the statement ``$\tilde X = X$, $\dbP$-a.s.''
we mean that these processes are equal to $dt\times d\dbP$-a.s.
When both of them are {c\`{a}dl\`{a}g}, clearly $\tilde X_t = X_t$,
$0\le t\le1$, $\dbP$-a.s.


\section{Second order target problem and relaxations}
\label{sect-2target}

In this section we start with the definitions and assumptions related
to the nonlinearity $H$ and its convex dual. Several spaces used in the
paper are also
introduced in Section~\ref{ss.def}. We then give the definition
of the original problem, two relaxed problems and the dual. We provide an
easy first string of inequalities in the final subsection.

\subsection{Definitions and assumptions}
\label{ss.def}

Let $H_t(\omega,y,z,\gamma)\dvtx  [0,1]\times\O\times\dbR\times\dbR
^d\times D_H
\to\dbR$ be $\dbF$-progressively measurable, where
$D_H\subset\dbR^{d\times d}$ is a given subset containing $0$.
We assume throughout the following:

\begin{assum}\label{assum-H}
For all $\omega\in\O$, $H$ is Lipschitz continuous
in $(y,z)$, uniformly in $(t,\omega, \gamma)$ and it is
uniformly continuous in $\omega$ under the $\dbL^\infty$-norm.
Moreover, we assume that it is
lower-semicontinuous in $\gamma$ and the conjugate $F$ defined at
(\ref{F}) below is measurable.
\end{assum}

In the sequel, we denote by $A\: B:=\operatorname{Tr}[A^{\mathrm{T}}B]$ for
$A, B\in\dbR^{d\times n}$. We introduce
the conjugate of $H$ with respect to $\gamma$ by
%
\begin{equation}\label{F}
F_t(\omega,y,z,a)
:=
\sup_{\gamma\in D_H} \biggl\{\frac{1}{2}a\dvtx \gamma
- H_t(\omega,y,z,\gamma)\biggr\}, \qquad a\in\dbS^{>0}_d.
\end{equation}
We notice that $F$ is measurable if $H$ is upper-semicontinuous (and
hence continuous) in $\gamma$ or if $D_H$ is compact; see, for
example,~\cite{BS}.
Moreover, since $H$ is uniformly continuous in $(\omega,y,z)$,
the domain of $F$ as a function of $a$ is independent of $(\omega, y,z)$.
Thus, we denote it by $D_{F_t}$.
By the uniform Lipschitz continuity of $H$ in $(y,z)$,
we know that
%
\begin{equation}
\label{FLip}
\begin{tabular}{p{290pt}@{}}
$F(\cd,a)$ is uniformly Lipschitz
continuous in $(y,z)$ and uniformly continuous in $\omega$,
uniformly on $(t,a)$, for every $a\in D_{F_t}$.
\end{tabular}
\end{equation}
Moreover, for our duality result of Section~\ref{sect-duality},
we need to further assume the following:
%
\begin{assum}\label{assum-F}
There is a constant $C$ such that,
for all $(t,\omega,y,z_1, z_2)$ and all $a\in D_{F_t}$:
\[
|F_t(\omega, y, z_1, a) -
F_t(\omega, y, z_2, a)|\le C|a^{1\slash2}(z_1-z_2)|.
\]
\end{assum}

We also define
%
\begin{equation}\label{Fhat}
\hat F_t(y,z)
:=
F_t(y,z,\hat a_t)\quad
\mbox{and}\quad
\hat F^0_t:=\hat F_t(0,0).
\end{equation}
In order to focus on our main idea, in this section we
shall restrict the probability measures in a subset
$\cP_H\subset\overline{\cP}_S$ defined below.
We will extend our results to more general cases,
as well as allowing $H$ to take value $\infty$,
in Section~\ref{sect-extension} below.

\begin{defn}
\label{defn-cP} Let $\cP_H$ denote the collection of
all those $\dbP\in\overline{\cP}_S$ such that
%
\begin{equation}\label{ellipticity}
\underline a_\dbP\le\hat a \le\overline a_\dbP,\qquad
dt\times d\dbP\mbox{-a.s. for some }
\underline a_\dbP, \overline a_\dbP\in\dbS^{>0}_d,
\end{equation}
and
\begin{equation}
\label{cP}
\dbE^\dbP\biggl[\int_0^1
(|\hat F^0_t|^2 + |H^0_t|^2)\,dt \biggr] <\infty.
\end{equation}
\end{defn}

\begin{rem}
\label{rem-cPHk}
In our accompanying paper~\cite{STZ09d} we consider a slightly
more general class $\cP_H^\k$ with a parameter $\k\in(1,2]$. The
$\cP_H$ in this paper coincides with the case $\k=2$ there. All the
results in this paper can be easily extended to the general case $\k
\in(1,2]$. In particular, Theorem~\ref{thm-duality} and Proposition
\ref{prop-sup=esssup} in this paper still hold true for general $\k$,
which are used in~\cite{STZ09d}, Theorem 4.6.
\end{rem}

It is clear that $\hat a_t\in D_{F_t}$,
$dt\times d\dbP$-a.s. for all $\dbP\in\cP_H$,
and by (\ref{FLip}) together with Assumption~\ref{assum-F},
%
\begin{eqnarray}\label{FLipschitz}
|\hat F_t(y_1, z_1)-\hat F_t(y_2,z_2)|\le C\bigl(|y_1-y_2|
+ |\hat a_t^{1\slash2}(z_1-z_2)|\bigr),
\nonumber
\\[-8pt]
\\[-8pt]
\eqntext{dt\times d\dbP\mbox{-a.s. for all } \dbP\in\cP_H.}
\end{eqnarray}

\begin{rem}
The Lipschitz continuity in $z$ in (\ref{FLipschitz}) is implied by the
following condition on $H$:
\begin{eqnarray*}
|H_t(y,z_1,\gamma)-H_t(y,z_2,\gamma)|
&\le&
C|\hat a_t^{1/2}(z_1-z_2)|,\qquad dt\times d\dbP\mbox{-a.s.}
\end{eqnarray*}
for some constant $C$ which does not depend on $(t,\omega,y,\gamma)$.
\end{rem}

We conclude this subsection by introducing
the spaces which will be needed for the formulation
of the second order target problems.
For any domain $D$ in an Euclidean space with appropriate
dimension,\vadjust{\goodbreak}
let $\dbL^0(D)$ denote the space of all
$\cF_1$-measurable random variables taking
values in $D$, and
$\dbH^0(D)$ the space of all $\dbF^+$-progressively
measurable processes taking values in $D$. Notice that here we use the
right limit filtration $\dbF^+$.
For any $\dbP\in\overline{\cP}_W$, let $\dbD^0(\dbP,D) $ be
the subspace of $\dbH^0(D)$ whose elements have
{c\`{a}dl\`{a}g} paths, $\dbP$-a.s.; $\dbI^0(\dbP,D)$ the subspace of
$\dbD^0(\dbP,D)$ whose elements $K$ have nondecreasing paths with
$K_0=0$, $\dbP$-a.s.; and $\dbS^0(\dbP,D)$ the subspace of
$\dbD^0(\dbP,D)$ whose elements have continuous paths, $\dbP$-a.s.

Moreover, let
%
\begin{eqnarray}\label{Hp}
\dbL^2(\dbP,D) &:=&
\{\xi\in\dbL^0(D)\dvtx  \dbE^{\dbP}[|\xi|^2]<\infty\},\nonumber\\
\dbH^2(\dbP,D)
&:=&
\biggl\{H\in\dbH^0(D)\dvtx  \dbE^{\dbP}\biggl[\int_0^1
|H_t|^2 \,dt\biggr]<\infty\biggr\},\nonumber\\
\dbD^2(\dbP,D)
&:=&
\Bigl\{Y\in\dbD^0(\dbP,D)\dvtx  \dbE^{\dbP}\Bigl[\sup_{0\le t\le1}
|Y_t|^2\Bigr]<\infty\Bigr\},\\
\dbI^2(\dbP,D)&:=& \dbD^2(\dbP,D)\cap\dbI^0(\dbP,D),\nonumber\\
\dbS^2(\dbP,D)& :=& \dbD^2(\dbP,D)\cap\dbS^0(\dbP,D),\nonumber
\end{eqnarray}
and denote
\begin{eqnarray*}
\hat\dbL^2_H(D) := \bigcap_{\dbP\in\cP_H}
\dbL^2(\dbP,D),\qquad \hat\dbH^2_H(D) :=
\bigcap_{\dbP\in\cP_H}\dbH^2(\dbP, D),
\end{eqnarray*}
and the corresponding subsets of {c\`{a}dl\`{a}g},
continuous processes, nondecreasing processes:
$\hat\dbD^2_H(D) := \bigcap_{\dbP\in\cP_H}\dbD^2(\dbP, D)$,
$\hat\dbS^2_H(D) := \bigcap_{\dbP\in\cP_H}\dbS^2(\dbP, D)$,\break
$\hat\dbI^2_H(D) := \bigcap_{\dbP\in\cP_H}\dbI^2(\dbP, D)$.

Finally, let
\begin{eqnarray*}
\hat\dbG^2_H(D_H) := \bigcap_{\dbP\in\cP_H} \dbG^2(\dbP
,D_H)\quad\mbox{and}\quad\hSM_H(\dbR^d) :=
\bigcap_{\dbP\in\cP_H}\SM_H(\dbP, \dbR^d),
\end{eqnarray*}
where
\[
\dbG^2(\dbP, D_H)
:=
\bigl\{\Gamma\in\dbH^0(D_H)\dvtx  \tfrac12\hat{a}\:\Gamma
-H(0,0,\Gamma)\in\dbH^2(\dbP,\dbR)
\bigr\}
\]
and $\cS\cM^2_H(\dbP,\dbR^d)\subset\dbD^2(\dbP,\dbR^d)$
is the space of all square integrable $(\dbP, \dbF^+)$-semimartingales
$Z$ with $\G\in\dbG^2(\dbP, D_H)$, where $\G$ is
defined by $d\langle Z,B\rangle_t = \G_t \dvtx d\langle B\rangle_t$,
$\dbP$-a.s.

\begin{rem}\label{rem-FF+}We emphasize that in the above spaces
we require the processes to be $\dbF^+$-progressively measurable. This
is important because the process $V^+$ in (\ref{V+}) is in general
$\dbF^+$-progressively measurable. See also Proposition \ref
{prop-V+=V} and the paragraph before it.

However, for fixed $\dbP\in\overline{\cP}_S$, it follows from the
Blumenthal zero--one law that $\dbE^\dbP[\xi|\cF_t]=\dbE^\dbP[\xi
|\cF^+_t]$, $\dbP$-a.s. for any $t\in[0,1]$\vadjust{\goodbreak} and $\dbP$-integrable
$\xi$. In particular, this shows that any $\cF_t^+$-measurable random
variable has an $\cF_t$-measurable $\dbP$-modification. Consequently,
for any fixed $\dbP$, we may view the processes in $\dbL^2(\dbP,D)$
as $\dbF$-progressively measurable.
\end{rem}

\subsection{The second order target problem}

For $Z\in\hSM_H(\dbR^d)$, it follows from
Karandikar~\cite{Karandikar} that It\^o's stochastic integrals
\[
\int_0^t Z_s \,dB_s\quad \mbox{and}\quad \int_0^t B_s \,dZ_s\qquad
\mbox{are defined } \cP_H\mbox{-q.s.}
\]
In particular, the quadratic covariation between $Z$ and $B$ is
well defined $\cP_H$-q.s. and has a density process $\G$:
%
\begin{equation}\label{ZG}
d\langle Z,B\rangle_t
=
\G_t d\langle B\rangle_t
= \G_t \hat a_t \,dt,\qquad \cP_H\mbox{-q.s.}
\end{equation}
For any $y\in\dbR$ and $Z\in\hSM_H(\dbR^d)$,
let $Y := Y^{y,Z}\in\hat\dbS^2_H(\dbR)$ denote
the controlled process defined by the following ODE (with random
coefficients):
%
\begin{eqnarray}\label{Yy}
Y_t &=& y- \int_0^t H_s(Y_s,Z_s,\G_s)\,ds
+\int_0^t Z_s\circ dB_s\nonumber\\
&=&y + \int_0^t \biggl(\frac{1}{2}\hat a_s\:\G_s
-H_s(Y_s,Z_s,\G_s)\biggr)\,ds\\
&&{} +\int_0^t Z_s\,dB_s,\qquad
t\in[0,1], \cP_H\mbox{-q.s.,}\nonumber
\end{eqnarray}
where $\circ$ denotes the Stratonovich stochastic integral.
We note that the well-posedness of (\ref{Yy}) follows
directly from the assumptions that
$\G\in\hat\dbG^2_H(D_H)$, $Z$ is square integrable
under each $\dbP\in\cP_H$, and $H$ is
uniformly Lipschitz continuous in $(y,z)$.

Let $\xi\in\dbL^0(\dbR)$.
Following Soner and Touzi~\cite{ST}, we introduce the second order
stochastic target problem:
%
\begin{eqnarray}\label{V}
\cV(\xi):= \inf\{y\dvtx  Y^{y, Z}_1 \ge\xi,
\cP_H\mbox{-q.s. for some } Z \in\hSM_H(\dbR^d)\}.
\end{eqnarray}

\subsection{Relaxations}
\label{sectrelaxation}

We relax the target problem (\ref{V}) by removing
the constraint that $Z$ is a semimartingale. For
any $y\in\dbR$, $\bar Z\in\hat\dbH^2_H(\dbR^d)$,
$\bar\G\in\hat\dbG^2_H(D_H)$, and $\dbP\in\cP_H$, let $\bar Y:=
\bar Y^{\dbP, y, \bar Z,\bar\G}\in\dbS^2(\dbP,\dbR)$ denote the
unique solution of
%
\begin{eqnarray}\label{Yybar}
\bar Y_t
=
y + \int_0^t \biggl(\frac{1}{2}\hat a_s\:\bar\G_s-H_s(\bar Y_s,
\bar Z_s,\bar\G_s)\biggr)\,ds
+\int_0^t\bar Z_s\,dB_s,
\nonumber
\\[-8pt]
\\[-8pt]
\eqntext{ t\in[0,1], \dbP\mbox{-a.s.}}
\end{eqnarray}
Here, we observe that the stochastic integral $\int_0^t Z_s\,dB_s$
may not have a $\cP_H$-q.s. version, in general,
and thus we can only define (\ref{Yybar}) under each $\dbP\in\cP_H$.

Our relaxed target problem is
%
\begin{eqnarray}
\label{Vbar}
\bar{\cV}(\xi)&:=& \inf\{y\dvtx  \exists(\bar Z,\bar\Gamma)\in
\hat\dbH^2_H(\dbR^d)\times\hat\dbG^2_H(D_H) \mbox{ such that}
\nonumber
\\[-8pt]
\\[-8pt]
\nonumber
&&\hspace*{57pt} \bar Y_1^{\dbP,y,
\bar Z,\bar\G} \ge\xi, \dbP\mbox{-a.s. for all } \dbP\in\cP
_H\}.
\end{eqnarray}

The main duality result of this paper relies\vspace*{-1.5pt} on the following
further relaxation of the above target problems.
For $y\in\dbR$, $\bar{\hspace*{-1pt}\bar Z}\in\hat\dbH^2_H(\dbR^d)$
and $\dbP\in\cP_H$, let $\bar{\hspace*{-1pt}\bar Y}:=\bar{\hspace*{-1pt}\bar Y}{}^{\dbP,
y,\bar{\hspace*{-1pt}\bar Z}}\in\dbS^2(\dbP,\dbR)$ be the unique solution of
%
\begin{equation}\label{Yyhat}
\bar{\hspace*{-1pt}\bar Y}_t = y + \int_0^t \hat F_s(\bar{\hspace*{-1pt}\bar Y}_s,
\bar{\hspace*{-1pt}\bar Z}_s)\,ds +\int_0^t\bar{\hspace*{-1pt}\bar Z}_s\,dB_s,\qquad
t\in[0,1], \dbP\mbox{-a.s.,}
\end{equation}
where existence and uniqueness of $\bar{\hspace*{-1pt}\bar Y}$
follows from (\ref{cP}) and (\ref{FLipschitz}). Here, again,
the stochastic integral $\int_0^t \bar{\hspace*{-1pt}\bar Z}_s\,dB_s$ may
not have a $\cP_H$-q.s. version. Our further relaxed
second order target problem does not involve the processes $\G$
and $\bar\G$, and is defined by
%
\begin{equation}
\label{Vhat}
\qquad\bar{\hspace*{-1pt}\bar\cV}(\xi):= \inf\{y\dvtx  \exists\bar{\hspace*{-1pt}\bar Z}
\in\hat\dbH^2_H(\dbR^d) \mbox{ s.t. } \bar{\hspace*{-1pt}\bar Y}{}^{\dbP, y,\bar{\hspace*{-1pt}\bar Z}}_1
\ge\xi, \dbP\mbox{-a.s. for all } \dbP\in\cP_H\}.
\end{equation}

\subsection{Dual formulation}
By (\ref{0-1 MRP}), each $\dbP\in\cP_H\subset\overline{\cP}_S$ satisfies
the martingale representation property. Let $\t$ be an $\dbF
$-stopping time and $\eta$ an $\cF_\t$-measurable and $\dbP$-square
integrable random variable. By (\ref{cP}),
(\ref{FLipschitz}) and the standard BSDE theory, the following BSDE
has a unique solution $(\cY^\dbP(\t,\eta), \cZ^\dbP(\t,\eta
))\in\dbS^2(\dbP,\dbR)\times
\dbH^2(\dbP,\dbR^d)$:
%
\begin{eqnarray}
\label{YZa}
\cY^\dbP_t(\t,\eta) &= &\eta- \int_t^\t\hat F_s(\cY^\dbP_s(\t
,\eta),\cZ^\dbP_s(\t,\eta))\,ds
\nonumber
\\[-8pt]
\\[-8pt]
\nonumber
&&{}-\int_t^\t\cZ^\dbP_s(\t,\eta) \,dB_s,
\qquad\dbP\mbox{-a.s.}
\end{eqnarray}
Now for any
$\xi\in\hat\dbL^2_H(\dbR)$, our dual formulation is
%
\begin{equation}
\label{vxi}
v(\xi) := \sup_{\dbP\in\cP_H} \cY^\dbP_0(1,\xi).
\end{equation}
By the Blumenthal zero--one law (\ref{0-1 MRP}), we know $\cY^\dbP
_0(1,\xi)$ is a constant, and thus~$v(\xi)$ is deterministic.

Our main focus of this paper is to provide conditions
which guarantee that the problems $\bar\cV(\xi)$,
$\bar{\hspace*{-1pt}\bar\cV}(\xi)$ and $v(\xi)$ agree. In order to
connect these problems to $\cV(\xi)$, we will need
an appropriate reformulation; see Section~\ref{sectweak2target}.

\subsection{Some preliminary results}
In this subsection we prove a straightforward string of
inequalities.
%
\begin{prop}
\label{p.easy}
Let Assumptions~\ref{assum-H} and~\ref{assum-F} hold true. Then, for
any $\xi\in\hat\dbL^2_H(\dbR)$,
%
\begin{equation}\label{VgebarVgehatV}
\cV(\xi) \ge \bar{\cV}(\xi) = \bar{\hspace*{-1pt}\bar\cV}(\xi)
\ge v(\xi).
\end{equation}
\end{prop}

\begin{pf}
(i) The first inequality holds true by definition of $\cV$ and
$\bar{\cV}$.\vspace*{-6pt}
\begin{longlist}[(iii)]
\item[(ii)] To prove that $\bar{\cV}(\xi)\ge\bar{\hspace*{-1pt}\bar\cV}(\xi)$,
let $y\in\dbR$, $\bar Z\in\hat\dbH^2_H(\dbR^d)$
and $\bar\G\in\hat\dbG^2_H(D_H)$ be such that
$\bar Y^{\dbP,y,\bar Z, \bar\G}_1 \ge\xi$,
$\dbP$-a.s. for all $\dbP\in\cP_H$. By the definition of the
conjugate function $F$,
\[
\tfrac{1}{2}\hat a_s\dvtx \bar\G_s-H_s(y, \bar Z_s,\bar\G_s) \le\hat
F_s(y, \bar Z_s)\qquad
\mbox{for all }
y\in\dbR.
\]
By the comparison theorem for ODEs, we
conclude that $\bar Y^{\dbP,y,\bar Z, \bar\G}_1
\le\bar{\hspace*{-1pt}\bar Y}{}^{\dbP, y,\bar Z}_1$, $\dbP$-a.s.
Thus, $\bar{\hspace*{-1pt}\bar Y}{}^{\dbP, y,\bar Z}_1\ge\xi, \dbP$-a.s. and,
therefore, $y\ge\bar{\hspace*{-1pt}\bar\cV}(\xi)$.

\item[(iii)] Similarly, to see that $\bar{\hspace*{-1pt}\bar\cV}(\xi)\ge\bar{\cV
}(\xi)$,
we consider some $y>\bar{\hspace*{-1pt}\bar\cV}(\xi)$ so that there exists
$\bar{\hspace*{-1pt}\bar Z}\in\hat\dbH^2_H(\dbR^d)$ such that
\[
\bar{\hspace*{-1pt}\bar Y}{}^{\dbP,y\bar{\hspace*{-1pt}\bar Z}}_1\ge\xi,
\qquad\dbP\mbox{-a.s. for all }
\dbP\in\cP_H.
\]
Then, for any $\eps>0$, it follows from the
lower-semicontinuity of $H$ in $\gamma$ that there exists a
progressively measurable process
$\bar\G\in\dbH^0(D_H)$ such that
\[
\hat F(\bar{\hspace*{-1pt}\bar Y},\bar{\hspace*{-1pt}\bar Z})-\eps
\le
\tfrac12\hat a\:\bar\G-H(\bar{\hspace*{-1pt}\bar Y},\bar{\hspace*{-1pt}\bar Z},\bar\G
)
\le
\hat F(\bar{\hspace*{-1pt}\bar Y},\bar{\hspace*{-1pt}\bar Z}).
\]
Then, $\bar\G\in\hat\dbG^2_H(D_H)$ and it follows from
classical estimates on ODEs that there exists a constant
$C$ such that, with $\bar y:=y+C\eps$, we have
\[
\bar Y^{\dbP,\bar y,\bar{\hspace*{-1pt}\bar Z},\bar\G}_1
\ge
\bar{\hspace*{-1pt}\bar Y}{}^{\dbP,y\bar{\hspace*{-1pt}\bar Z}}_1\ge\xi
\qquad\mbox{a.s. for all }
\dbP\in\cP_H.
\]
Hence, $\bar y\ge\bar{\cV}(\xi)$.
Since $\varepsilon>0$ and $y>\bar{\hspace*{-1pt}\bar\cV}(\xi)$ are arbitrary,
we conclude that $\bar{\hspace*{-1pt}\bar\cV}(\xi)\ge{\bar\cV}(\xi)$.

\item[(iv)] The final inequality $\bar{\hspace*{-1pt}\bar\cV}(\xi)\ge v(\xi)$ can
be proved similarly to (ii) above by using the comparison theorem for
BSDEs.\quad\qed
\end{longlist}
\noqed\end{pf}

\begin{rem}
Consider the Markovian case $H_t(y,z,\gamma)=
h(t,B_t,y,z,\gamma)$ and $\xi=g(B_1)$, for some deterministic
functions $h,g$.
Assume in addition that the PDE (\ref{pdeu}) has a solution
$u\in C^{1,2}$ with appropriate growth. Then, by the classical
verification argument of stochastic control,
one can prove that $u(0,0)=v(\xi)$. Moreover, if $H$ is convex,
then it follows from a direct application of It\^o's formula that
$u(0,0)= \bar{\cV}(\xi) = \bar{\hspace*{-1pt}\bar\cV}(\xi) = v(\xi)$.
If in addition $\{Du(t,B_t),t\in[0,1]\}\in\hSM_H(\dbR^d)$,
then we also have $u(0,0)=\cV(\xi) = \bar{\cV}(\xi) = \bar{\hspace*{-1pt}\bar\cV}(\xi) = v(\xi)$.
Finally, any optimal $\dbP^*$ (if exists) for the problem $v(\xi)$ satisfies
\begin{eqnarray*}
\tfrac{1}{2}\hat{a}_t\dvtx  D^2u(t,B_t)-H_\cd(\cd, u, Du, D^2
u)(t,B_t)
=
F_\cd(\cd, u, Du, \hat{a}_\cd)(t,B_t),\qquad \dbP^*\mbox{-a.s.}
\end{eqnarray*}
\end{rem}

In the non-Markovian case, we shall prove in the
next section our main duality result
$\bar{\cV}(\xi)=\bar{\hspace*{-1pt}\bar\cV}(\xi) = v(\xi)$
and that the optimal $(\bar Z, \bar\G)$, $\bar{\hspace*{-1pt}\bar Z}$,
for the problems $\bar{\cV}(\xi)$ and $\bar{\hspace*{-1pt}\bar\cV}(\xi)$,
respectively, exist. However, we are not able to prove
${\cV}(\xi) = \bar{\cV}(\xi)$ in general. In order to obtain
a result of this type, we shall introduce a slight modification
of these problems by restricting $\dbP$ to smaller sets;
see Section~\ref{sectweak2target} below.

\section{The main results}
\label{sect-duality}

This section is devoted to the proof of reverse inequalities.

\subsection{Conditional expectation}
We first establish a dynamic programming principle to
prove our duality result $\bar{\hspace*{-1pt}\bar\cV}(\xi) = v(\xi)$.
The understanding of the regular conditional probability distributions
(r.c.p.d.)
is crucial for this result. Indeed, let $\dbP$ be an arbitrary probability
measure on $\O$ and $\t$ be an $\dbF$-stopping time.
By Stroock and Varadhan~\cite{SV},
there exists a r.c.p.d. $\dbP^\omega_\t$ for all $\omega\in\O$
satisfying the following:
\begin{itemize}[--]
\item[--] For each $\omega\in\O$, $\dbP^\omega_\t$ is a probability
measure on $\cF_1$;

\item[--] For each $E\in\cF_1$, the mapping $\omega\to\dbP^\omega_\t
(E)$ is
$\cF_\t$-measurable;

\item[--] For $\dbP$-a.e. $\omega\in\O$, $\dbP^\omega_\t$ is the
conditional probability measure of $\dbP$ on $\cF_\t$, that is,
for every bounded $\cF_1$-measurable random variable $\xi$ we have
\[
\dbE^\dbP(\xi|\cF_\t)(\omega) = \dbE^{\dbP^\omega_\t}(\xi),
\qquad\dbP\mbox{-a.s.};
\]
\item[--] For each $\omega\in\O$,
%
\begin{equation}
\label{Oto}
\dbP^\omega_\t(\O^{\omega}_\t) = 1 \qquad\mbox{where } \O^{\omega
}_\t
\:= \{\omega'\in\O\dvtx  \omega'(s) = \omega(s), 0\le s\le\t
(\omega)\}.
\end{equation}
\end{itemize}

The goal of this subsection is to understand
$\dbP^\omega_\t$ for $\dbP\in\cP_H$. Roughly, we shall prove that
$\dbP^\omega_\t$ satisfies the properties of Definition~\ref{defn-cP}
on a shifted space; see Lemma~\ref{lem-rcpdH} below.
To do that, we introduce some notation:
\begin{itemize}[--]
\item[--] For $0\le t\le1$, denote by $\O^t:= \{\omega\in C([t,1], \dbR
^d)\dvtx  \omega(t)=0\}$
the shifted canonical space; $B^{t}$ the shifted canonical process on
$\O^t$; $\dbP^{t}_0$ the shifted Wiener measure; $\dbF^{t}$ the
shifted filtration generated by $B^{t}$.

\item[--] For $0\le s\le t\le1$ and $\omega\in\O^s$, define the shifted
path $\omega^t\in\O^t$:
\[
\omega^t_r := \omega_r-\omega_t\qquad
\mbox{for all }
r\in[t, 1];
\]
\item[--] For $0\le s \le t\le1$ and $\omega\in\O^s$, $\tilde\omega\in
\O^t$,
define the concatenation path $\omega\otimes_{t} \tilde\omega\in
\O^s$ by
\begin{eqnarray*}
(\omega\otimes_t \tilde\omega) (r) := \omega_r\1_{[s,t)}(r) +
(\omega_{t} + \tilde\omega_r)\1_{[t, 1]}(r)
\qquad\mbox{for all }
r\in[s,1].
\end{eqnarray*}
\item[--] For $0\le s\le t \le1$ and an $\cF^{s}_{1}$-measurable
random variable $\xi$ on $\O^s$, for each $\omega\in\O^s$,
define the shifted $\cF^{t}_{1}$-measurable random variable $\xi^{t,
\omega}$ on $\O^t$ by
\[
\xi^{t, \omega}(\tilde\omega) :=\xi(\omega\otimes_t \tilde
\omega)\qquad
\mbox{for all }
\tilde\omega\in\O^t.
\]
Similarly, for an $\dbF^{s}$-progressively measurable
process $X$ on $[s, 1]$ and $(t,\omega)\in[s,1]\times\O^s$,
the shifted process $\{X^{t, \omega}_r, r\in[t,1]\}$ is $\dbF
^{t}$-progressively measurable.

\item[--] For $\dbF$-stopping time $\t$, we shall simplify the notation
as follows:
\[
\omega\otimes_\t\tilde\omega:= \omega\otimes_{\t(\omega)}
\tilde\omega, \qquad\xi^{\t,\omega} := \xi^{\t(\omega),\omega
}, \qquad X^{\t,\omega} := X^{\t(\omega),\omega}.
\]
\end{itemize}

The r.c.p.d. $\dbP^\omega_\t$ induces naturally a probability
measure $\dbP^{\t,\omega}$ on $\cF^{\t(\omega)}_1$ such that the
$\dbP^{\t,\omega}$-distribution of $B^{\t(\omega)}$ is equal to the
$\dbP^\omega_\t$-distribution of $\{B_t-B_{\t(\omega)}, t\in[\t
(\omega),1]\}$.
By (\ref{Oto}), it is clear that
for every bounded and $\cF_1$-measurable random variable $\xi$,
%
\begin{equation}
\label{Pto}
\dbE^{\dbP^\omega_\t}[\xi] = \dbE^{\dbP^{\t,\omega}}[\xi
^{\t,\omega}] .
\end{equation}
We shall also call $\dbP^{\t,\omega}$ the r.c.p.d. of $\dbP$.

For $0\le t\le1$, following the same arguments as in
Section~\ref{sect-setup} but restricting to the canonical space
$\O^t$, we may define martingale measures $\dbP^{t, \a}$
for each $\dbF^{t}$-progressively measurable
$\dbS^{>0}_d$-valued process $\a$ such that
$\int_{t}^{1}|\a_r| \,dr <\infty$, $\dbP^{t}_0$-a.s. Let
$\overline\cP_S^{t}$ denote the set of all such measures
$\dbP^{t, \a}$. Similarly, we may define the density process
$\hat a^{t}$ of the quadratic variation process $\langle B^{t}\rangle$.

We first have the following result.

\begin{lem}
\label{lem-rcpd}
Let $\dbP\in\overline{\cP}_S$ and $\t$ be an $\dbF$-stopping
time. Then, for
$\dbP$-a.e. $\omega\in\O$, $\dbP^{\t,\omega} \in\overline
{\cP}^{\t(\omega)}_S$ and
%
\begin{equation}
\label{arcpd}
\qquad\hat a_s^{\t,\omega}(\tilde\omega) = \hat a^{\t(\omega
)}_s(\tilde\omega) \qquad\mbox{for } ds \times d\dbP^{\t,\omega
}\mbox{-a.e.}\ (s, \tilde\omega)\in[\t(\omega),1]\times\O^{\t
(\omega)},
\end{equation}
where the left-hand side above is the shifted process of original density
process $\hat a$ on $\O=\O_0$ and the right-hand side is the density
process on the shifted space~$\O^{\t(\omega)}$.
\end{lem}

\begin{pf} The proof of $\dbP^{\t,\omega} \in\overline{\cP}^{\t
(\omega)}_S$ is relegated to the \hyperref[sectappendix]{Appendix}. We now prove (\ref{arcpd}).

Since $d\langle B_\cd- B_{\t}\rangle_t = \hat a_t \,dt$, $\dbP$-a.s., then
$d\langle B_\cd- B_{\t}\rangle_t = \hat a_t \,dt$, $\dbP^\omega_\t$-a.s.
for $\dbP$-a.e. $\omega\in\O$.
Note that, for each $\omega\in\O$ and $t\ge\t(\omega)$,
\[
\hat a_t(\omega) = \hat a_t\bigl(\omega\otimes_{\t} \omega^{\t
(\omega)}\bigr) = \hat a^{\t, \omega}_t\bigl(\omega^{\t(\omega)}\bigr).
\]
This implies that $d\langle B^{\t(\omega)}_\cd\rangle_t = \hat
a^{\t
,\omega}_t \,dt$, $\dbP^{\t,\omega}$-a.s. for $\dbP$-a.e. $\omega
\in\O$. Now (\ref{arcpd}) follows from the definition of $\hat a^{\t
(\omega)}$.
\end{pf}

We next study the r.c.p.d. for $\dbP\in\cP_H$. For each $(t,\omega
)\in[0,1]\times\O$, let
%
\begin{eqnarray}
\label{HFrcpd}
H^{t,\omega}_s(\tilde\omega, y,z, \gamma) &:=& H_s(\omega\otimes
_t\tilde\omega, y,z,\gamma),
\nonumber
\\[-8pt]
\\[-8pt]
\nonumber
 \hat F^{t,\omega}_s(\tilde
\omega, y,z) &:= &F_s\bigl(\omega\otimes_t\tilde\omega, y, z, \hat
a^{t}_s(\tilde\omega)\bigr)
\end{eqnarray}
for all $(s, \tilde\omega)\in[t,1]\times\O^t$ and $(y,z,\gamma
)\in\dbR\times\dbR^d\times D_H$. We emphasize that in the
definition of $\hat F^{t,\omega}$ we use the density process $\hat
a^t$ in the shifted space.\vadjust{\goodbreak} This is important in (\ref{hatFcont}) below.
However, by Lemma~\ref{lem-rcpd} we actually have
\begin{eqnarray}
\hat F^{t,\omega}_s(\tilde\omega, y,z) = F_s\bigl(\omega\otimes
_t\tilde\omega, y, z, \hat a^{t,\omega}_s(\tilde\omega)\bigr) = \hat
F_s(\omega\otimes_t\tilde\omega, y, z),\nonumber\\
\eqntext{ds \times d\dbP^{t,\omega}\mbox{-a.e. } (s, \tilde\omega)\in
[t,1]\times\O^t, \dbP\mbox{-a.e. } \omega\in\O.}
\end{eqnarray}
Since $H$ and $F$ are uniformly continuous in $\omega$ under the
$\dbL^\infty$-norm, by Assumption~\ref{assum-H} and (\ref{FLip}), we
also have
%
\begin{eqnarray}
\label{hatFcont}
\begin{tabular}{p{270pt}@{}}
$H^{t,\omega}_s(\tilde\omega, y,z, \gamma)$ and $\hat
F^{t,\omega}_s(\tilde\omega, y,z)$
are uniformly continuous in $\omega$ under the $\dbL
^\infty$-norm.
\end{tabular}
\end{eqnarray}
We remark that $F_s(\omega\otimes_t\tilde\omega, y, z, \hat
a^{t,\omega}_s(\tilde\omega))$ is in general not continuous in
$\omega$ because~$\hat a$ is not continuous in $\omega$, in
general; see Lemma~\ref{lem-atilde}. Similarly, as a consequence of
(\ref{hatFcont}), we see that for any $\dbP^t\in\overline{\cP}^t_S$,
%
\begin{eqnarray}
\label{someiffall}
\dbE^{\dbP^t}\biggl[\int_t^1 \bigl(|H^{t,\omega}_s(0)|^2 + |\hat
F^{t,\omega}_s(0)|^2\bigr)\,ds\biggr]<\infty
\nonumber
\\[-8pt]
\\[-8pt]
\eqntext{\mbox{for some $\omega\in\O$ iff it holds for all }
\omega\in\O.}
\end{eqnarray}
We now extend Definition~\ref{defn-cP} to the shifted space.
%
\begin{defn}
\label{defn-cPt} Let $\cP^t_H$ denote the collection of all those
$\dbP\in\overline{\cP}^t_S$ such that
%
\begin{eqnarray}\label{cPt}
&&\underline a_\dbP\le\hat a^t \le\overline a_\dbP,\ ds\times
d\dbP\mbox{-a.e. on } [t,1]\times\O^t \mbox{ for some } \underline
a_\dbP, \overline a_\dbP\in\dbS^{>0}_d,
\nonumber
\\
&&\hspace*{61pt}\dbE^\dbP\biggl[\int_t^1 \bigl(|H^{t,\omega}_s(0)|^2 +
|\hat F^{t,\omega}_s(0)|^2\bigr)\,ds\biggr]<\infty\\
\eqntext{\mbox{for all or,
equivalently, some } \omega\in\O.}
\end{eqnarray}
\end{defn}

Then we have the following.
%
\begin{lem}
\label{lem-rcpdH}
Let Assumption~\ref{assum-H} hold true. Then, for any $\dbF$-stopping
time $\t$ and $\dbP\in\cP_H$, the r.c.p.d. $\dbP^{\t,\omega}
\in\cP_H^{\t(\omega)}$, for $\dbP$-a.e. $\omega\in\O$.
\end{lem}

\begin{pf} Let $\dbP= \dbP^\a\in\cP_H\subset\overline{\cP}_S$. By
Lemma~\ref{lem-rcpd} we have $\dbP^{\t,\omega} \in\overline{\cP
}^{\t(\omega)}_S$, $\dbP$-a.s. By (\ref{ellipticity}) and (\ref
{cP}), it holds for $\dbP$-a.e. $\omega\in\O$ that
\begin{eqnarray*}
&\displaystyle\underline{a}_\dbP\le\hat a_s^{\t,\omega}(\tilde\omega) \le
\overline{a}_\dbP, ds\times d\dbP^{\t,\omega}\mbox{-a.e. } (s,\tilde\omega)\in[\t(\omega),1]\times\O^{\t(\omega
)},&\\
&\displaystyle \dbE^{\dbP^{\t,\omega}}\biggl[\int_{\t(\omega)}^1
\bigl(\bigl|F_s\bigl(\omega\otimes_{\t} \tilde\omega, 0, 0, \hat a^{\t,\omega
}_s(\tilde\omega)\bigr)\bigr|^2 + |H_s(\omega\otimes_{\t}\tilde\omega,
0,0,0)|^2\bigr)\,ds\biggr]<\infty.&
\end{eqnarray*}
This, together with (\ref{arcpd}) and (\ref{HFrcpd}), implies (\ref
{cPt}), and thus completes the proof.
\end{pf}

We remark that in this paper we actually use the r.c.p.d. only on
deterministic times. However, the r.c.p.d. on stopping times will be
important in our accompanying paper~\cite{STZ09d}.

\subsection{The duality result}

To establish our main duality result, we need
the following assumption on the terminal data.

\begin{assum}
\label{assum-xi}
$\xi$ is uniformly continuous in $\omega$ under the $\dbL^\infty$-norm.
\end{assum}

Under Assumptions~\ref{assum-H} and~\ref{assum-xi}, there exists a
modulus of continuity function~$\rho$ for $\xi$ and $H$ in $\omega
$. Then, for any $0\le t\le s\le1$, $(y,z)\in[0,1]\times\dbR\times
\dbR^d$, and $\omega, \omega'\in\O$, $\tilde\omega\in\O^t$,
\begin{eqnarray*}
|\xi^{t,\omega}(\tilde\omega) - \xi^{t,\omega'}(\tilde
\omega)| &\le&\rho(\|\omega-\omega'\|_t) \quad\mbox{and}\\
|\hat F^{t,\omega}_s(\tilde\omega, y,z) - \hat F^{t,\omega
'}_s(\tilde\omega, y,z)| &\le&\rho(\|\omega-\omega'\|_t),
\end{eqnarray*}
where $\|\omega\|_t:= \sup_{0\le s\le t}|\omega_s|$, $0\le t\le1$.
We next define for all $\omega\in\O$,
%
\begin{eqnarray}
\label{Lambda}
\L(\omega) := \sup_{0\le t\le1}\L_t(\omega)
\nonumber
\\[-8pt]
\\[-8pt]
\eqntext{\mbox{where } \L
_t(\omega):= \displaystyle\mathop{\sup}_{\dbP\in\cP^t_H} \biggl(\dbE^\dbP\biggl[|\xi
^{t,\omega}|^2 + \int_t^1 |\hat F^{t,\omega}_s(0)|^2\,ds\biggr]
\biggr)^{1/2} .}
\end{eqnarray}
By (\ref{hatFcont}) and following the same arguments as for (\ref
{someiffall}), we have
%
\begin{equation}
\label{Lambda-finite}
\L(\omega) <\infty \mbox{ for some } \omega\in\O\mbox{ iff it
holds for all } \omega\in\O.
\end{equation}
Moreover, when $\L$ is finite, it is uniformly continuous in $\omega
$ under the $\dbL^\infty$-norm and is therefore $\cF_1$-measurable.

Our main duality result is as follows:

\begin{thmm}
\label{thm-duality}
Let Assumptions~\ref{assum-H},~\ref{assum-F},~\ref{assum-xi} hold,
and assume further that
%
\begin{equation}
\label{Lambda-integrability}
\dbE^\dbP[|\L|^2]<\infty\qquad \mbox{for all } \dbP\in\cP_H.
\end{equation}
Then ${\bar\cV}(\xi)=\bar{\hspace*{-1pt}\bar\cV}(\xi) = v(\xi)$, and
existence holds for the problem $\bar{\hspace*{-1pt}\bar\cV}(\xi)$. Moreover, if~$F$ has a progressively\vspace*{1pt} measurable optimizer,
existence also holds for the problem $\bar\cV(\xi)$.
\end{thmm}

We first provide several examples that satisfy the
hypothesis of the theorem and then prove it in Section~\ref{ss.proof}.

\subsection{Examples}
\label{sect-examples}

\begin{example}[(Linear generator)]
Assume that $H$ is linear in $\gamma$:
\[
H_t(y,z,\gamma)
=
f_t(y,z)+\tfrac12\sigma_t\sigma_t^{\mathrm{T}}\:\gamma,\vadjust{\goodbreak}
\]
where $f_t(y,z)$ and $\sigma_t$ satisfy appropriate conditions for our
assumptions to hold. Notice that the domain of $F$ is reduced to a
one-point set:
\[
F_t(y,z,a)
=
f_t(y,z)\1_{\{a=\sigma_t\sigma_t^{\mathrm{T}}\}}+\infty\1_{\{a\neq\sigma
_t\sigma_t^{\mathrm{T}}\}}.
\]
Then, the present formulation of the second order target problem is
clearly equivalent to the classical formulation under the reference
measure $\dbP^{\sigma\sigma^{\mathrm{T}}}$ which ignores any uncertainty on
the diffusion coefficient.
\end{example}

\begin{example}[(Uncertain volatility models)]
\label{ex Peng}
Set $H_t(y,z,\gamma) := G(\gamma) := \frac{1}{2}[\bar\si^2
\gamma^+ - \underline{\si}^2 \gamma^-]$, where $\bar\si
>\underline{\si}\ge0$. This is the context studied by Denis and
Martini~\cite{DM}. By straightforward calculation, we find $\Dom
(F_t)= [\underline{\si}^2, \bar\si^2]$, and for any $a\in
[\underline{\si}^2, \bar\si^2]$, $F(a) = 0$. It is easily seen that
all our assumptions are satisfied. Moreover, we have $\bar{\cV}(\xi)
= \bar{\hspace*{-1pt}\bar{\cV}}(\xi) = \dbE^G(\xi)$ for appropriate
random variable $\xi$, where $\dbE^G$ is the $G$-expectation defined
in Peng~\cite{Peng-G}. More connections between this paper and
$G$-martingales are established in our accompanying paper~\cite{STZ09b}.
\end{example}

\begin{example}[(Hedging under gamma constraints)]
\label{ex gammaconstraint}
Let $\underline{\Gamma},\overline{\Gamma}\ge0$ be two given constants.
The problem of superhedging under Gamma constraint, as introduced in
\cite{STgammaSIAM,CST} and~\cite{STgamma}, corresponds to
the specification $H_s(y,z,\gamma)=H(\gamma)=\frac12\sigma^2\gamma
$ for $\gamma\in[-\underline{\Gamma},\overline{\Gamma}]$, and
$+\infty$ otherwise. By straightforward calculation, we see that
$F(a)=\frac12(\overline{\Gamma}(a-\sigma^2)^++\underline
{\Gamma}(a-\sigma^2)^-)$. If both bounds are finite, the
domain of the dual function $F$ is the nonnegative real line. The dual
formulation of this paper coincides with that of~\cite{STgamma}.
\end{example}

\subsection{Proof of the duality result}
\label{ss.proof}

The rest of this section is devoted to the proof of Theorem \ref
{thm-duality}. From now on, we shall always assume Assumptions~\ref
{assum-H}, \ref{assum-F},~\ref{assum-xi} and
that (\ref{Lambda-integrability}) hold. In particular, we notice that
(\ref{Lambda-integrability}) and (\ref{Lambda-finite}) imply that
%
\begin{equation}
\label{Lambda-finite2}
\L_t(\omega) <\infty\qquad \mbox{for all } (t,\omega)\in[0,1]\times
\O.
\end{equation}

To prove the theorem, we define the following value process $V_t$ pathwise:
%
\begin{equation}\label{Vtomega}
V_t(\omega)
:=
\sup_{\dbP\in\cP^t_H} \cY^{\dbP,t,\omega}_t(1,\xi)\qquad \mbox
{for all } (t,\omega)\in[0,1]\times\O,
\end{equation}
where, for any $(t_1, \omega)\in[0,1]\times\O$, $\dbP\in\cP
_H^{t_1}$, $t_2\in[t_1,1]$, and any $\eta\in\dbL^2(\dbP,\cF
_{t_2})$, we denote $\cY^{\dbP,t_1,\omega}_{t_1}(t_2,\eta
):= y^{\dbP,t_1,\omega}_{t_1}$, where $(y^{\dbP,t_1,\omega
},z^{\dbP,t_1,\omega})$ is the solution to the following BSDE on the
shifted space $\O^{t_1}$ under $\dbP$:
%
\begin{eqnarray}
\label{YPt}
y^{\dbP,t_1,\omega}_s = \eta^{t_1,\omega} - \int_s^{t_2}
\hat F^{t_1,\omega}_r(y^{\dbP,t_1,\omega}_r, z^{\dbP
,t_1,\omega}_r)
\,dr - \int_s^{t_2} z^{\dbP,t_1,\omega}_r \,dB^{t_1}_r,
\nonumber
\\[-8pt]
\\[-8pt]
\eqntext{d s\in
[t_1,t_2], \dbP\mbox{-a.s.}}
\end{eqnarray}
In view of the Blumenthal zero--one law (\ref{0-1 MRP}), $\cY^{\dbP
,t,\omega}_t(1,\xi)$ is
constant for any given $(t, \omega)$ and $\dbP\in\cP^t_H$. Moreover,
since $\omega_0=0$ for all $\omega\in\O$, it is clear that, for
the $\cY^\dbP$ defined in (\ref{YZa}),
\[
\cY^{\dbP,0,\omega}(t,\eta) = \cY^{\dbP}(t,\eta) \quad\mbox{and}\quad
V_0(\omega) = v(\xi) \qquad\mbox{for all } \omega\in\O.
\]

\begin{lem}
\label{lem-hatYreg}
Assume all the conditions in Theorem~\ref{thm-duality} hold. Then for
all $(t,\omega)\in[0,1]\times\O$, we have $|V_t(\omega)|\le C\L
_t(\omega)$. Moreover, for all $(t,\omega,\omega')\in[0,1]\times
\O^2$, $|V_t(\omega)-V_t(\omega')|\le C\rho(\|\omega-\omega'\|
_t)$. Consequently, $V_t$ is $\cF_t$-measurable for every $t\in[0,1]$.
\end{lem}

\begin{pf}(i) For each $(t,\omega)\in[0,1]\times\O$ and $\dbP\in\cP
^t_H$, on $[t,1]$ we have
\begin{eqnarray}
y^{\dbP, t,\omega}_s = \xi^{t,\omega} - \int_s^1 [\hat
F^{t,\omega}_r(0) + \gamma_s y^{\dbP, t,\omega}_r + z^{\dbP,
t,\omega}_r (\hat a^t_r)^{1/2} \eta_r^T] \,dr - \int_s^1
z^{\dbP, t,\omega}_r \,dB^t_r,\nonumber\\
\eqntext{\dbP\mbox{-a.s.},}
\end{eqnarray}
where $\gamma, \eta$ are bounded, thanks to (\ref{FLipschitz}). Define
%
\begin{equation}
\label{M}
M_s := \exp\biggl(-\int_t^s \eta_r \,dB^t_r - \int_t^s \biggl[\gamma_r +
{1\over2}| (\hat a^t_r)^{1/2}\eta_r^T|^2\biggr]\,dr\biggr).
\end{equation}
Applying It\^{o}'s formula, we obtain
\begin{eqnarray*}
y^{\dbP, t,\omega}_t = M_ty^{\dbP, t,\omega}_t = M_T \xi
^{t,\omega} - \int_t^1 M_s \hat F^{t,\omega}_s(0) \,ds -\int_t[\cds
]\,dB^t_s,\qquad \dbP\mbox{-a.s.}
\end{eqnarray*}
Thus,
\begin{eqnarray*}
|y^{\dbP, t,\omega}_t|^2 &=& \biggl|\dbE^\dbP\biggl[M_T \xi
^{t,\omega} - \int_t^1 M_s \hat F^{t,\omega}_s(0) \,ds\biggr]
\biggr|^2\\
&\le&\biggl|\dbE^\dbP\biggl[\sup_{t\le s\le T}M_s |\xi^{t,\omega
}| + \int_t^1|\hat F^{t,\omega}_s(0)| \,ds\biggr]\biggr|^2\\
&\le& C\dbE^\dbP\Bigl[\sup_{t\le s\le T}|M_s|^2\Bigr] \dbE^\dbP
\biggl[ |\xi^{t,\omega}|^2 + \int_t^1|\hat F^{t,\omega}_s(0)|^2
\,ds\biggr].
\end{eqnarray*}
Since $\gamma, \eta$ are bounded, by standard arguments we see that
\[
|y^{\dbP,t,\omega}_{t}|^2\le C\dbE^\dbP\biggl[|\xi^{t,\omega}|^2
+ \int_t^1 |\hat F^{t,\omega}_s(0)|^2\,ds\biggr] \le C|\L_t(\omega)|^2.
\]
Since $\dbP\in\cP^t_H$ is arbitrary, we get $|V_t(\omega)|\le C\L
_t(\omega)$.

(ii) Similarly, for $(t,\omega,\omega')\in[0,1]\times\O^2$ and
$\dbP\in\cP^t_H$, denote
\begin{eqnarray*}
\delta y& :=& y^{\dbP,t,\omega}-y^{\dbP,t,\omega'},\qquad
\delta z := z^{\dbP,t,\omega}_t-z^{\dbP,t,\omega'},\qquad
\delta\xi:= \xi^{t,\omega}-\xi^{t,\omega'},\\
\delta F &:=& \hat F^{t,\omega}-\hat F^{t,\omega'}.
\end{eqnarray*}
Then, for $ s\in[t,1]$, $|\delta\xi| + |\delta F_s|\le C\rho(\|
\omega-\omega'\|_t)$ and
\begin{eqnarray}
\delta y_s
=
\delta\xi
- \int_s^1 [\delta F_r (y^{\dbP,t,\omega}_r, z^{\dbP,t,\omega
}_r)+ \tilde\gamma_r \delta y_r + \delta z_r(\hat a^t_r)^{1/2}\tilde\eta_r^T] \,dr
- \int_s^1 \delta z_r \,dB^t_r,\nonumber\\
\eqntext{\dbP\mbox{-a.s.},}
\end{eqnarray}
where $\tilde\gamma, \tilde\eta$ are bounded, thanks to (\ref
{FLipschitz}) again. Define $\tilde M$ as in (\ref{M}) but corresponding
to $(\tilde\gamma, \tilde\eta)$. Then following the arguments in
(i), we obtain $|\delta y_t|\le C\rho(\|\omega-\omega'\|_t)$.
Since $\dbP$ is arbitrary, we prove the lemma.
\end{pf}

The following dynamic programming principle plays
a central role in our analysis.

\begin{prop}\label{prop-DPP}
Assume all the conditions in Theorem~\ref{thm-duality} hold. Then
\begin{eqnarray*}
V_{t_1}(\omega)
=
\sup_{\dbP\in\cP^{t_1}_H} \cY^{\dbP,t_1,\omega}_{t_1}
(t_2,V^{t_1,\omega}_{t_2})\qquad
\mbox{for all } 0\le t_1<t_2\le1 \mbox{ and } \omega\in\O.
\end{eqnarray*}
\end{prop}

\begin{pf} To simplify the presentation, we assume without loss of
generality that $t_1=0$ and $t_2=t$. That is, we shall prove
%
\begin{equation}
\label{DPP}
v(\xi) = \sup_{\dbP\in\cP_H} \cY^{\dbP}_0(t,V_{t}).
\end{equation}
Denote $(y^\dbP, z^\dbP) := (\cY^\dbP(1,\xi), \cZ^\dbP(1,\xi))$.
\begin{longlist}[(ii)]
\item[(i)] For any $\dbP\in\cP_H$, note that
\[
y^{\dbP}_s
=
y^{\dbP}_{t} - \int_s^{t} \hat F_r(y^{\dbP}_r, z^{\dbP}_r) \,dr -
\int_s^{t} z^{\dbP}_r \,dB_r, \qquad s\in[0, t], \dbP\mbox{-a.s.}
\]
By Lemma~\ref{lem-rcpdH}, for $\dbP$-a.e. $\omega\in\O$, the
r.c.p.d. $\dbP^{t,\omega}\in\cP_H^{t}$. Since solutions of BSDEs
can be constructed via Picard iteration, one can easily check that
%
\begin{equation}
\label{rcpdBSDE}
y^{\dbP}_{t}(\omega) = \cY^{\dbP^{t,\omega}, t, \omega
}_{t}(1,\xi)\qquad \mbox{for } \dbP\mbox{-a.e. } \omega\in\O.
\end{equation}
Then by the definition of $V_{t}$ we get
%
\begin{equation}
\label{YP<V}
y^{\dbP}_{t}(\omega) \le V_{t}(\omega) \qquad \mbox{for } \dbP\mbox{-a.e. } \omega\in\O.
\end{equation}
It follows from the comparison principle for BSDEs that $y^{\dbP}_{0}
\le\cY^{\dbP}_{0}(t,V_t)$. Since $\dbP\in\cP_H$ is arbitrary,
this shows that $v(\xi)\le\sup_{\dbP\in\cP_H}
\cY^{\dbP}_{0}(t,V_t)$.\vspace*{1pt}

\item[(ii)] It remains to prove the other inequality. Fix $\dbP\in\cP
_H$ and arbitrary $\e>0$. Since $\O$ is separable, there exists a
partition $E_t^i\in\cF_t, i=1,2,\ldots$ such that $\|\omega-\omega
'\|_t\le\e$ for any $i$ and any $\omega,\omega'\in E_t^i$. For
each $i$, fix an $\hat\omega_i\in E_t^i$, and let $\dbP^i_t\in\cP
^t_H$ be an $\e$-optimizer of $V_t(\hat\omega_i)$, that is, $
V_t(\hat\omega_i) \le\cY^{\dbP^i_t, t, \hat\omega_i}_t + \e$.

For each $n\ge1$, define $\dbP^{n}:= \dbP^{n,\e}$ by
%
\begin{eqnarray}
\label{Pn}
\dbP^{n}(E) := \dbE^{\dbP}\Biggl[\sum_{i=1}^n \dbE^{\dbP^i_t}[(\1
_E)^{t,\omega}]\1_{E^i_t}\Biggr] + \dbP(E\cap\hat E^n_t)
\nonumber
\\[-8pt]
\\[-8pt]
\eqntext{\mbox{where } \hat E^n_t \dfnn\displaystyle\bigcup_{i>n} E^i_t.}
\end{eqnarray}
That is, $\dbP^{n}=\dbP$ on $\cF_t$, and its r.c.p.d. $(\dbP
^n)^{t,\omega} = \dbP^i_t$ for $\omega\in E^i_t$, $1\le i\le n$,
and $(\dbP^n)^{t,\omega} = \dbP^{t,\omega}$ for $\omega\in\hat
E^n_t$. We claim that
%
\begin{equation}
\label{PncPS}
\dbP^n\in\cP_H.
\end{equation}
The proof is similar to Lemmas~\ref{lem-rcpd} and~\ref{lem-rcpdH},
and thus is also postponed to the \hyperref[sectappendix]{Appendix}.

Now for $1\le i\le n$ and $\omega\in E^i_t$, by Lemma \ref
{lem-hatYreg} and its proof we see that
\begin{eqnarray*}
V_t(\omega) &\le& V_t(\hat\omega_i) + C\rho(\e) \le\cY^{\dbP
^i_t, t, \hat\omega_i}_t(1,\xi) + \e+ C\rho(\e)\\
&\le& \cY^{\dbP^i_t, t, \omega}_t(1,\xi) + \eps+ C\rho(\e) =
\cY^{(\dbP^{n})^{t,\omega}, t, \omega}_t(1,\xi) + \e+ C\rho(\e).
\end{eqnarray*}
Here as usual the constant $C$ varies from line to line. Then it
follows from (\ref{rcpdBSDE}) that
%
\begin{eqnarray}
\label{V<YP}
V_t \le y^{\dbP^{n}}_t + \e+ C\rho(\e), \qquad \dbP^{n}\mbox{-a.s.
on } \bigcup_{i=1}^n E^i_t.
\end{eqnarray}
Let $(y^n, z^n):=(y^{n,\e}, z^{n,\e})$ denote the solution to the
following BSDE on $[0,t]$:
\begin{eqnarray}
y^{n}_s = [y^{\dbP^{n}}_t + \e+ C\rho(\e)]\1_{\bigcup
_{i=1}^n E^i_t} + V_t\1_{\hat E^n_t} - \int_s^t \hat F_r(y^{n}_r,
z^{n}_r)\,dr - \int_s^tz^{n}_r \,dB_r,\nonumber\\
\eqntext{\dbP\mbox{-a.s}.}
\end{eqnarray}
By the comparison principle of BSDEs we know $\cY^{\dbP}_0(t,V_t) \le
y^{n}_0$. Since $\dbP^n = \dbP$ on~$\cF_t$, we have
\[
y^{\dbP^n}_s = y^{\dbP^n}_t - \int_s^t \hat F_r(y^{\dbP^n}_r,
z^{\dbP^n}_r)\,dr - \int_s^tz^{\dbP^n}_r \,dB_r,\qquad  s\in[0,t], \dbP
\mbox{-a.s.}
\]
By the standard arguments in BSDE theory we get
\[
|y^n_0 - y^{\dbP^n}_0|^2 \le C\dbE^\dbP[|\e+ C\rho(\e)|^2 +
|V_t - y^{\dbP^n}_t|^2 \1_{\hat E^n_t}].
\]
By Lemma~\ref{lem-hatYreg} and its proof we have $|V_t|\le C\L_t$ and
$|y^{\dbP^n}_t|\le C\L_t$, $\dbP$-a.s.
Then
\begin{eqnarray*}
\cY^{\dbP}_0(t,V_t) &\le& y^{n}_0
\le
y^{\dbP^n}_0 + C\bigl(\e+ \rho(\e)\bigr) + C(\dbE^\dbP[|\L_t|^2 \1
_{\hat E^n_t}])^{1/2}\\
&\le&
v(\xi) + C\bigl(\e+ \rho(\e)\bigr) + C(\dbE^\dbP[|\L_t|^2 \1_{\hat
E^n_t}])^{1/2}.
\end{eqnarray*}
Recall (\ref{Lambda-integrability}) and notice that $\hat
E^n_t\downarrow\varnothing$. By sending $n\to\infty$ and applying the
dominated convergence theorem we get
\[
\cY^{\dbP}_0(t,V_t) \le v(\xi)+ C\bigl(\e+ \rho(\e)\bigr)\qquad \mbox{for
all } \dbP\in\cP_H.
\]
Since $\e>0$ is arbitrary, we complete the proof.\quad\qed
\end{longlist}
\noqed\end{pf}

We next introduce the right limit of the $V$ which is defined for each
$(t,\omega)$ and is clearly $\dbF^+$-progressively measurable:
\[
V^+_t := \mathop{\limsup}_{r\in\dbQ\cap(t,1], r\downarrow t}V_r.
\]

\begin{lem}
\label{lem-V+}
Assume all the conditions in Theorem~\ref{thm-duality} hold. Then
%
\begin{equation}
\label{V+}
\qquad V^+_t = \lim_{r\in\dbQ\cap(t,1], r\downarrow t}V_r,\qquad \cP_H\mbox{-q.s. and, thus, } V^+ \mbox{ is {{c\`{a}dl\`{a}g}} } \cP_H\mbox{-q.s.}
\end{equation}
\end{lem}

\begin{pf} For each $\dbP\in\cP_H$, denote
\[
\tilde V^\dbP:= V - \cY^\dbP(1,\xi).
\]
Then $\tilde V^\dbP_t \ge0$, $\dbP$-a.s. For any $0\le t_1 < t_2\le
1$, let $(y^{\dbP, t_2}, z^{\dbP, t_2}) := (\cY^{\dbP}(t_2,
V_{t_2}), \break \cZ^{\dbP}(t_2, V_{t_2}))$. Note that $\cY^{\dbP
}_{t_1}(t_2, V_{t_2})(\omega) = \cY_{t_1}^{\dbP, t_1, \omega
}(t_2, V^{t_1,\omega}_{t_2})$ for $\dbP$-a.s. $\omega$.
Then by Proposition~\ref{prop-DPP} we get $V_{t_1} \ge y^{\dbP,
t_2}_{t_1}$, $\dbP$-a.s. Notice that $y^{\dbP,1} = y^\dbP$. Denote
\[
\tilde y^{\dbP, t_2}_t := y^{\dbP, t_2}_t - y^{\dbP}_t, \qquad\tilde
z^{\dbP, t_2}_t := \hat a_t^{-1\slash2}(z^{\dbP, t_2}_t - z^{\dbP}_t).
\]
Then $\tilde V^\dbP_{t_1} \ge\tilde y^{\dbP, t_2}_{t_1}$, $\dbP
$-a.s. and $(\tilde y^{\dbP, t_2}, \tilde z^{\dbP, t_2})$ satisfies
the following BSDE on $[0,t_2]$:
\[
\tilde y^{\dbP, t_2}_t = \tilde V^\dbP_{t_2} - \int_t^{t_2} f^\dbP
_s(\tilde y^{\dbP, t_2}_s, \tilde z^{\dbP, t_2}_s)\,ds - \int_t^{t_2}
\tilde z^{\dbP, t_2}_s \,dW^\dbP_s, \qquad\dbP\mbox{-a.s.},
\]
where
\begin{eqnarray*}
f^\dbP_t(\omega,y,z) &:=& \hat F_t\bigl(\omega, y+y^\dbP
_t(\omega), \hat a_t^{-1\slash2}(\omega)\bigl(z+z^\dbP_t(\omega
)\bigr)\bigr)\\
&&{}- \hat F_t(\omega, y^\dbP_t(\omega), \hat a_t^{-1\slash
2}(\omega)z^\dbP_t(\omega)).
\end{eqnarray*}
Notice that $f^{\dbP}_t(0,0) =0$, and $f^\dbP$ is uniformly Lipschitz
continuous in $(y,z)$. Following the definition in~\cite{Peng} and
\cite{CP}, $\tilde V^\dbP$ is a weak $f^\dbP$-supermartingale under
$\dbP$. Now applying the downcrossing inequality Theorem 6 of \cite
{CP}, one can easily see that, for $\dbP$-a.e. $\omega$, the limit
$\lim_{r\in\dbQ\cap(t,1],r\downarrow t} \tilde V^\dbP_r(\omega)$
exists for all $t\in[0,1]$. Note that~$y^\dbP$ is continuous, $\dbP
$-a.s. We get that the $\limsup$ in the definition of $V^+$ is in fact
the $\lim$, $\dbP$-a.s. Then,
\begin{eqnarray*}
V^+_t = \lim_{r\in\dbQ\cap(t,1], r\downarrow t}V_r, t\in[0,1]
\qquad \mbox{and, therefore, }
V^+ \mbox{ is {c\`{a}dl\`{a}g}}, \cP_H\mbox{-q.s.}
\end{eqnarray*}
\upqed\end{pf}

We are now ready to prove our main duality result.\vadjust{\goodbreak}

\begin{pf*}{Proof of Theorem \protect\ref{thm-duality}}
We proceed in several steps.

\textit{Step} 1.
We first show that $V^+$ is a strong $\hat F$-supermartingale under
each $\dbP\in\cP_H$. For any $\dbP\in\cP_H$, denote $\tilde
V^{+,\dbP} := V^+ - y^\dbP$. Given $0\le t_1 < t_2< 1$, let $r^1_n\in
\dbQ\cap(t_1,t_2], r^1_n \downarrow t_1$ and $r^2_n\in\dbQ\cap
(t_2,1], r^2_n \downarrow t_2$. We have $\tilde V^\dbP_{r^1_n} \ge
\tilde y^{\dbP, r^2_m}_{r^1_n}$, $\dbP$-a.s. for any $m, n\ge1$.
Sending $n\to\infty$, we get $\tilde V^{+,\dbP}_{t_1} \ge\tilde
y^{\dbP, r^2_m}_{t_1}$, $\dbP$-a.s. for any $m\ge1$. Sending $m\to
\infty$, by the stability of BSDEs we get $\tilde V^{+,\dbP}_{t_1}
\ge\tilde y^{+,\dbP, t_2}_{t_1}$, $\dbP$-a.s. where
\begin{eqnarray*}
\tilde y^{+,\dbP, t_2}_t = \tilde V^{+,\dbP}_{t_2} - \int_t^{t_2}
f^\dbP_s(\tilde y^{+,\dbP, t_2}_s, \tilde z^{+,\dbP, t_2}_s)\,ds -
\int_t^{t_2} \tilde z^{+,\dbP, t_2}_s \,dW^\dbP_s, \qquad\dbP\mbox{-a.s.}
\end{eqnarray*}
That is, $\tilde V^{+,\dbP}$ is also a weak $f^\dbP$-supermartingale
under $\dbP$. Applying Theorem~7 of~\cite{CP}, $\tilde V^{+,\dbP}$
is a strong $f^\dbP$-supermartingale under $\dbP$. That is, recalling
(\ref{0-1 MRP}), for any $\overline{\dbF}^\dbP$-stopping times $\t
_1, \t_2$ with $\t_1\le\t_2$, we have $\tilde V^{+,\dbP}_{\t_1}
\ge\tilde y^{+,\dbP, \t_2}_{\t_1}$, $\dbP$-a.s. where
\begin{eqnarray}
\tilde y^{+,\dbP, \t_2}_t = \tilde V^{+,\dbP}_{\t_2} - \int_t^{\t
_2} f^\dbP_s(\tilde y^{+,\dbP, \t_2}_s, \tilde z^{+,\dbP, \t_2}_s)\,ds
- \int_t^{\t_2} \tilde z^{+,\dbP, \t_2}_s \,dW^\dbP_s,\nonumber\\
\eqntext{t\in[0,\t
_2], \dbP\mbox{-a.s.}}
\end{eqnarray}
This implies that $V^+_{\t_1}\ge y^{+,\dbP,\t_2}_{\t_1}$, $\dbP
$-a.s. where $y^{+,\dbP,\t_2}_{t}:= \tilde y^{+,\dbP, \t_2}_t +
y^\dbP_t, z^{+,\dbP,\t_2}_{t}:= \hat a_t^{1\slash2}(\tilde
z^{+,\dbP, \t_2}_t + z^\dbP_t)$ satisfy
\begin{eqnarray*}
y^{+,\dbP, \t_2}_t = V^{+}_{\t_2} - \int_t^{\t_2} \hat F_s(\tilde
y^{+,\dbP, \t_2}_s, \tilde z^{+,\dbP, \t_2}_s)\,ds - \int_t^{t_2}
\tilde z^{+,\dbP, \t_2}_s \,d B_s, \qquad\dbP\mbox{-a.s.}
\end{eqnarray*}
That is, $V^+$ is a strong $\hat F$-supermartingale under $\dbP$.

\textit{Step} 2. For each $\dbP\in\cP_H$, applying the nonlinear
Doob--Meyer decomposition in~\cite{Peng}, there exist unique ($\dbP
$-a.s.) processes $\bar{\hspace*{-1pt}\bar Z} {}^\dbP\in\dbH^2(\dbP, \dbR^d)$ and
$K^\dbP\in\dbI^2(\dbP,\dbR)$ such that
%
\begin{eqnarray}
\label{V+strong}
V^+_t = V^+_0 + \int_0^t \hat F_s(V^+_s, \bar{\hspace*{-1pt}\bar Z} {}^\dbP_s)\,ds
+\int_0^t \bar{\hspace*{-1pt}\bar Z} {}^\dbP_s \,dB_s -K^\dbP_t,
\nonumber
\\[-8pt]
\\[-8pt]
\eqntext{0\le t\le1,
\dbP\mbox{-a.s.}}
\end{eqnarray}
Remark~\ref{rem-LX} below provides a simpler argument for this result.
By Karandi\-kar~\cite{Karandikar}, since $V^+$ is a {{c\`{a}dl\`{a}g}}
semimartingale under each $\dbP\in\cP_H$, we can define uniquely a
universal process $\bar{\hspace*{-1pt}\bar Z}$ by $d\langle V^+, B\rangle_t = \bar{\hspace*{-1pt}\bar
Z}_t d\langle B\rangle_t$, so that $\bar{\hspace*{-1pt}\bar Z} = \bar{\hspace*{-1pt}\bar Z} {}^\dbP$,
$dt\times d\dbP$-a.s. for all $\dbP\in\cP_H$. Thus, we have
%
\begin{eqnarray}\label{decomposition V}
V^+_t = V^+_0 + \int_0^t \hat F_s(V^+_s, \bar{\hspace*{-1pt}\bar Z}_s)\,ds +\int_0^t
\bar{\hspace*{-1pt}\bar Z}_s \,dB_s -K^\dbP_t,
\nonumber
\\[-8pt]
\\[-8pt]
\eqntext{0\le t\le1, \dbP\mbox{-a.s. for all } \dbP\in\cP_H.}
\end{eqnarray}

\textit{Step} 3. We remark that $V^+_0$ is $\cF^+_0$-measurable and
is not a constant in general. For each $\dbP\in\cP_H\subset
\overline\cP_S$, and each $r\in\dbQ\cap(0,1]$, we have $V_0 \ge
y^{\dbP, r}_0$, where $y^{\dbP, r}_0$ is a constant, thanks to the
Blumenthal zero--one law (\ref{0-1 MRP}) under $\dbP$. It is clear that
$\lim_{r\downarrow0}y^{\dbP, r}_0 = V^+_0$, $\dbP$-a.s. Then $V_0
\ge V^+_0$, $\dbP$-a.s. for all $\dbP\in\cP_H$.\vspace*{-1pt} Now by the
comparison of ODE and recalling (\ref{Yyhat}) and (\ref{decomposition
V}), we see that $\bar{\hspace*{-1pt}\bar Y}{}^{V_0, \bar{\hspace*{-1pt}\bar Z}}_1 \ge\bar{\hspace*{-1pt}\bar Y}{}^{V^+_0,
\bar{\hspace*{-1pt}\bar Z}}_1 \ge V^+_1 = \xi$, $\dbP$-a.s. for all
$\dbP\in\cP_H$.
Now by the definition of $\hat\cV(\xi)$, we get $\bar{\hspace*{-1pt}\bar{\cV
}}(\xi) \le V_0 = v(\xi)$. This, together with (\ref{VgebarVgehatV}),
proves $\bar{\cV}(\xi) = \bar{\hspace*{-1pt}\bar{\cV}}(\xi) = v(\xi)$.
Moreover, the process $\bar{\hspace*{-1pt}\bar Z}$ in (\ref{decomposition V}) is clearly
the optimal control for the problem $\bar{\hspace*{-1pt}\bar{\cV}}(\xi)$.
Finally, when $F$ has a progressively measurable optimizer, the
existence of the optimal control for the problem $\bar{\cV}(\xi)$ is obvious.
\end{pf*}

\begin{rem}
\label{rem-LX}
Following a suggestion of Nicole El Karoui, we derive the
decomposition (\ref{V+strong}) by the following alternative argument.
Consider the following reflected BSDE:
\begin{eqnarray*}
\cases{
\displaystyle\bar{\hspace*{-1pt}\bar Y} {}^\dbP_t = \xi-\int_t^1 \hat F_s(\bar{\hspace*{-1pt}\bar Y}{}^\dbP_s,
\bar{\hspace*{-1pt}\bar Z} {}^\dbP_s)\,ds
-\int_t^1 \bar{\hspace*{-1pt}\bar Z} {}^\dbP_s \,dB_s
+K^\dbP_1-K^\dbP_t,\vspace*{2pt}\cr
\bar{\hspace*{-1pt}\bar Y} {}^\dbP_t\ge V^+_t, [\bar{\hspace*{-1pt}\bar Y} {}^\dbP_{t-} - V^+_{t-}]\,d
K^\dbP_t = 0.}\qquad
0\le t\le1,
\dbP\mbox{-a.s.}
\end{eqnarray*}
By Lepeltier and Xu~\cite{LX}, the above RBSDE has a unique solution
and $\bar{\hspace*{-1pt}\bar Y} {}^\dbP$ is {c\`{a}dl\`{a}g}. Then it suffices to show
that $\bar{\hspace*{-1pt}\bar Y} {}^\dbP= V^+$, $\dbP$-a.s. In fact, if they are not equal,
without loss of generality we assume $\bar{\hspace*{-1pt}\bar Y} {}^\dbP_0 > V^+_0$.
For each $\e>0$, denote $\t_\e:= \inf\{t\dvtx  \bar{\hspace*{-1pt}\bar Y}{}^\dbP_t
\le V^+_t+\e\}$. Then $\t_\e$ is an $\overline{\dbF}^\dbP
$-stopping time and $\bar{\hspace*{-1pt}\bar Y}{}^{\dbP}_{t-}\ge V^+_{t-}+\e>
V^+_{t-}$ for all $t\le\t_\e$. Then $K^\dbP_t = 0$, $t\le\t_\e$,
and thus
\[
\bar{\hspace*{-1pt}\bar Y} {}^\dbP_t = \bar{\hspace*{-1pt}\bar Y} {}^\dbP_{\t_\e} -\int_t^{\t_\e
} \hat F_s(\bar{\hspace*{-1pt}\bar Y} {}^\dbP_s, \bar{\hspace*{-1pt}\bar Z} {}^\dbP_s)\,ds
-\int_t^{\t_\e} \bar{\hspace*{-1pt}\bar Z} {}^\dbP_s \,dB_s.
\]
Note that $\bar{\hspace*{-1pt}\bar Y} {}^\dbP_{\t_\e} \le V^+_{\t_\e}+\e$, by
comparison theorem for BSDEs and following standard arguments we have
$\bar{\hspace*{-1pt}\bar Y} {}^\dbP_0 \le y^{+,\dbP,\t_\e}_0 + C\e\le V^+_0 + C\e
$. Since $\e$ is arbitrary, this contradicts with $\bar{\hspace*{-1pt}\bar Y} {}^\dbP
_0 > V^+_0$.
\end{rem}

We conclude this section by establishing a representation formula for
$V^+$, which will be important for our accompanying paper \cite
{STZ09d}. For each $\dbP\in\cP_H$ and $t\in[0,1]$, denote
%
\begin{eqnarray}
\label{cPtP}
\cP_H(t, \dbP)& :=& \{\dbP'\in\cP_H\dvtx  \dbP' = \dbP\mbox{ on } \cF
_t\}\quad \mbox{and}
\nonumber
\\[-8pt]
\\[-8pt]
\nonumber
 \quad\cP_H(t+, \dbP) &:= &\{\dbP'\in\cP_H\dvtx  \dbP' =
\dbP\mbox{ on } \cF^+_t\}.
\end{eqnarray}
Then we have the following.
%
\begin{prop}\label{prop-sup=esssup}
Assume all the conditions in Theorem~\ref{thm-duality} hold. Then, for
each $\dbP\in\cP_H$,
\begin{eqnarray*}
V_t& =& {\esup\limits_{\dbP' \in\cP_H(t, \dbP)}} ^{\hspace*{-5pt} \dbP} \cY
^{\dbP'}_t(1,\xi) \quad\mbox{and}\\
V^+_t &= &\mathop{\esup}_{\dbP' \in\cP
_H(t+, \dbP)} {}^{\hspace*{-10pt} \dbP} \cY^{\dbP'}_t(1,\xi),\qquad\dbP
\mbox{-a.s.}
\end{eqnarray*}
\end{prop}

\begin{pf} Fix $\dbP\in\cP_H$. Denote
\begin{eqnarray*}
V^\dbP_t := {\esup\limits_{\dbP' \in\cP_H(t, \dbP)}} ^{ \hspace*{-5pt}\dbP}
\cY^{\dbP'}_t(1,\xi)\quad \mbox{and}\quad V^{\dbP,+}_t := {\esup\limits_{\dbP'
\in\cP_H(t+, \dbP)}} ^{\hspace*{-8pt} \dbP} \cY^{\dbP'}_t(1,\xi).
\end{eqnarray*}
\vspace*{-20pt}
\begin{longlist}[(ii)]
\item[(i)] We first prove the equality for $V$. For each $\dbP'\in\cP_H(t,
\dbP)\subset\cP_H$, by (\ref{YP<V}) we have $y^{\dbP'}_t \le V_t$,
$\dbP'$-a.s. Since $\dbP' = \dbP$ on $\cF_t$, then $y^{\dbP'}_t
\le V_t$, $\dbP$-a.s. and, thus, $V^\dbP_t \le V_t$, $\dbP$-a.s.

On the other hand, proceeding as in step (ii) of the proof of
Proposition~\ref{prop-DPP}, we define $\dbP^n$ for each $n, \e$ by
(\ref{Pn}). By (\ref{PncPS}), it is clear that $\dbP^n \in\cP
_H(t,\dbP)$. Then it follows from (\ref{V<YP}) that
\begin{eqnarray}
\dbP[V_t \le V^{\dbP}_t + \e+ C\rho(\e)] \ge\dbP
[V_t \le y^{\dbP^n}_t + \e+ C\rho(\e)]\ge\dbP\biggl[\bigcup_{1\le
i\le n}E^i_t\biggr] \to1 \nonumber\\
\eqntext{\mbox{as } n\to\infty.}
\end{eqnarray}
That is, $V_t \le V^{\dbP}_t + \e+ C\rho(\e)$, $\dbP$-a.s. for
all $\e>0$. This implies that $V_t \le V^{\dbP}_t$, $\dbP$-a.s.

\item[(ii)] We\vspace*{1pt} now prove the equality for $V^+$. First, for each $\dbP'\in
\cP_H(t+, \dbP)\subset\cP_H$ and $r\in\dbQ\,\cap\,(t,1]$, we have
$y^{\dbP'}_r \le V_r$, $\dbP'$-a.s. Sending $r\downarrow t$, we
obtain $y^{\dbP'}_t \le V^+_t$, $\dbP'$-a.s. Since both $y^{\dbP
'}_t$ and $V^+_t$ are $\cF^+_t$-measurable and $\dbP'=\dbP$ on $\cF
^+_t$, then $y^{\dbP'}_t \le V^+_t$, $\dbP$-a.s. and, thus, $V^{\dbP
,+}_t \le V_t$, $\dbP$-a.s.

On the other hand, for each $r\in\dbQ\cap(t,1]$, since $V_r = V^\dbP
_r$, $\mathbb{P}$-a.s. Following the same arguments in~\cite{STZ09d} Theorem
4.3, Step (iii) (we emphasize that there is no danger of cycle proof
here!), we have
%
\begin{eqnarray}
\label{Pn-Neveu}
\mbox{there exist } \dbP_n\in\cP(r, \dbP) \mbox{ such that } \cY
^{\dbP_n}_r(1,\xi) \uparrow V_r, \dbP\mbox{-a.s.}
\end{eqnarray}
Then, it follows from the stability of BSDEs that
\[
\cY^\dbP_t(r, V_r) = \cY^\dbP_t\Bigl(r, \lim_{n\to\infty}\cY
^{\dbP_n}_r(1,\xi)\Bigr) = \lim_{n\to\infty}\cY^\dbP_t(r,
\cY^{\dbP_n}_r(1,\xi)).
\]
Since $\dbP_n \in\cP(r, \dbP) \subset\cP(t+, \dbP)$, we have
\[
\cY^\dbP_t(r, V_r) = \lim_{n\to\infty}\cY^{\dbP_n}_t(r, \cY
^{\dbP_n}_r(1,\xi)) = \lim_{n\to\infty}\cY^{\dbP_n}_t(1,\xi
) \le V^{\dbP, +}_t,\qquad \dbP\mbox{-a.s.}
\]
Sending $r\downarrow t$, by the stability of BSDEs again we obtain
$V^+_t \le V^{\dbP, +}_t$, $\dbP$-a.s.\quad\qed
\end{longlist}
\noqed\end{pf}

After the completion of this paper, Marcel Nutz provides us the
following result which shows that, under our conditions that $F$ and
$\xi$ are uniformly continuous in~$\omega$, actually $V^+=V$.
However, we decide to keep our original\vadjust{\goodbreak} arguments because they are
applicable to more general cases, for example, the case in Section \ref
{sectweak2target}
where we do not require the uniform continuity of $F$ and $\xi$.

\begin{prop}[(M. Nutz)]
\label{prop-V+=V}
Assume all the conditions in Theorem~\ref{thm-duality} hold. Then
$V^+_t = V_t$, $\cP_H$-q.s.
\end{prop}

\begin{pf}First, by Lemma~\ref{lem-hatYreg} $V^+$ is uniformly continuous
in $\omega$ with the same modulus of continuity function $\rho$.
Since $V^+$ is $\dbF^+$-progressively measurable, for any $\d>0$, we
have $|V^+_t(\omega)-V^+_t(\omega')|\le C\rho(\|\omega-\omega
'\|_{t+\d})$. Sending $\d\to0$, we get $|V^+_t(\omega
)-V^+_t(\omega')|\le C\rho(\|\omega-\omega'\|_{t})$ and, thus,
$V^+$ is $\dbF$-progressively measurable.

By Proposition~\ref{prop-sup=esssup}, it is clear that $V^+_t \le
V_t$, $\cP_H$-q.s. On the other hand, for any $\dbP\in\cP_H$ and
$\dbP'\in\cP_H(t,\dbP)$, by the second equality of Proposition \ref
{prop-sup=esssup} we have $\cY^{\dbP'}_t(1,\xi) \le V^+_t$, $\dbP
'$-a.s. Since both sides of above are $\cF_t$-measurable and $\dbP' =
\dbP$ on $\cF_t$, we have $\cY^{\dbP'}_t(1,\xi) \le V^+_t$, $\dbP
$-a.s. Then the first equality of Proposition~\ref{prop-sup=esssup}
implies that $V_t \le V^+_t$, $\dbP$-a.s. Therefore, $V^+=V$, $\cP_H$-q.s.
\end{pf}

\section{A weaker version of the second order target problem}
\label{sectweak2target}

The purpose of this section is to suggest a slight modification of the
second order stochastic target problem so that its value is not
affected by the relaxations of Section~\ref{sectrelaxation}. The key
tool for this is the aggregation approach developed in our accompanying
paper~\cite{STZ09a}. The idea is to restrict our attention to an
(uncountable) subset of $\cP_H$, constructed out of a countable
subset, so that a dominating measure is available.

As a consequence of this modified setup, we shall remove the continuity
assumption on $\xi$. However, we still assume the nonlinearity $H$
satisfies Assumption~\ref{assum-H}, in particular, $H$ is uniformly
continuous in $\omega$ and the domain $D_{F_t}$ of its convex
conjugate $F$ is deterministic; see Section~\ref{sect-extension} for
the general case.

\subsection{\texorpdfstring{The dominating probability measure $\hP$}{The dominating probability measure P}}
\label{secthatP}

Throughout this section we fix a \textit{countable} subset $T_0\subset
[0,1]$ containing the end-points $\{0,1\}$, together with a \textit{countable}
sequence $\cA_0:=\{\a^i, i\ge1\}$ of deterministic
integrable mappings $\a^i \dvtx  [0,1] \to\dbS^{>0}_d$ satisfying the
concatenation property:
%
\begin{equation}\label{concatenation}
\a^i\1_{[0,t)}+\a^j\1_{[t,1]}\in\cA_0\qquad \mbox{for all } i, j\ge1,
t\in T_0.
\end{equation}
Note that $\a^i$ is deterministic, then by Lemma~\ref{lem-atilde},
$\hat a = \a^i$, $\dbP^{\a^i}$-a.s., and, thus, $\cA_0$ is a
generating class of diffusion coefficients in the sense of Definition~4.7
in~\cite{STZ09a}. Following Definition 4.8 in~\cite{STZ09a}, let
$\cA$ be the separable class of diffusion coefficients generated by
$(\cA_0, T_0)$. Following Proposition 8.3 in~\cite{STZ09a}, let $\cP
(\cA) \subset\overline\cP_S$ denote the corresponding measures.\vadjust{\goodbreak}
Then, by Definition 4.8 in~\cite{STZ09a},
%
\begin{eqnarray}
\label{PA}\qquad
\dbP\in\cP(\cA) \quad\mbox{if and only if}\quad \hat a = \sum
_{n=0}^\infty\sum_{i=1}^\infty\a^i \1_{E^n_i} \1_{[\t_n, \t
_{n+1})},\qquad \dbP\mbox{-a.s.}
\end{eqnarray}
for some
\begin{itemize}
\item sequence of $\dbF$-stopping times $\{\t_n, n\ge0\}$ with
values in $T_0$, with $\tau_0=0$, $\t_n<\t_{n+1}$ on $\{\t_n<1\}$,
and $\inf\{n\dvtx  \t_n =1\} <\infty$,

\item and some partition $\{E^n_i,i\ge1\}\subset\cF_{\t_n}$ of $\O$.
\end{itemize}
Finally, we assume $\dbP^i:= \dbP^{\a^i}\in\cP_H$ and denote $\cP
^\cA_H := \cP(\cA)\cap\cP_H$.

The dominating measure is now defined by
%
\begin{equation}\label{Phat}
\hP:= \hP^{\cA_0,T_0} := \sum_{i=1}^\infty2^{-i}\dbP^{i}.
\end{equation}
Clearly, $\hP$ is a dominating measure of $\{\dbP^i, i\ge1\}$.
By Proposition 4.11 in~\cite{STZ09a}, $\hP$ is in fact a dominating
measure of $\cP(\cA)$, and thus of $\cP^\cA_H$.
Therefore, $\cP^\cA_H$-q.s. reduces to $\hP$-a.s.

\subsection{\texorpdfstring{The second order target problem under $\hP$}{The second order target problem under P}}\label
{sect2targethatP}

Recall the spaces defined in~(\ref{Hp}). Let $\hat\dbL^2_0(D) :=
\bigcap_{i\ge1} \dbL^2(\dbP^i,D)$, and define the spaces $\hat\dbH
^2_0(D)$, $\hat\dbD^2_0(D)$, $\hat\dbS^2_0(D)$, $\hat\dbI
^2_0(D)$, $\hat\dbG^2_0(D_H)$ and $\hSM_0(\dbR^d)$ similarly.

Now for an $\cF_1$-measurable r.v. $\xi$, the modified second order
target problem under $\hP$ is
%
\begin{equation}\label{V0}
{\cV}_0(\xi):= \inf\{y\in\dbR\dvtx  Y^{y,Z}_1 \ge\xi, \hP
\mbox{-a.s. for some } Z\in\hSM_0(\dbR^d)\},
\end{equation}
where $Y^{y,Z}\in\hat\dbS^2_0(\dbR)$ is defined by (\ref{Yy}),
except that $\cP_H$-q.s. is replaced with $\hP$-a.s. (or,
equivalently, $\cP^\cA_H$-q.s.).

Next, notice that the families of processes $\{\bar Y^{\dbP^i}, i\ge
1\}$ and $\{\bar{\hspace*{-1pt}\bar Y}{}^{\dbP^i}, i\ge1\}$, defined by (\ref
{Yybar}) and (\ref{Yyhat}) respectively, can be aggregated into
processes $\bar Y$ and $\bar{\hspace*{-1pt}\bar Y}$, thanks to Theorem 5.1 in \cite
{STZ09a}. We then define the following relaxations of (\ref{V0}):
%
\begin{eqnarray}\label{Vbar0}
\bar{\cV}_0(\xi)
&:=&
\inf\{y\dvtx  \bar Y_1^{y,\bar Z,\bar\G} \ge\xi, \hP\mbox{-a.s. for some } (\bar Z,\bar\G)\in\hat\dbH^2_0(\dbR^d)\times
\hat\dbG^2_0(D_H)\},\hspace*{-35pt}\\
\label{Vhat0}
\bar{\hspace*{-1pt}\bar\cV}_0(\xi)
&:=&
\inf\{y\dvtx  \bar{\hspace*{-1pt}\bar Y}{}^{y,\bar{\hspace*{-1pt}\bar Z}}_1 \ge\xi, \hP
\mbox{-a.s. for some } \bar{\hspace*{-1pt}\bar Z}\in\hat\dbH^2_0(\dbR^d)
\}.
\end{eqnarray}
%
Finally, our modified dual formulation under $\hP$ is
%
\begin{equation}\label{v0}
v_0(\xi) := \sup_{\dbP\in\cP^\cA_H} \cY^\dbP_0(1,\xi),
\end{equation}
where $\cY^\dbP$ is defined by means of the BSDE (\ref{YZa}).
Similar to (\ref{VgebarVgehatV}), it is obvious that
%
\begin{equation}\label{VgebarVgehatV0}
\cV_0(\xi) \ge \bar\cV_0(\xi) = \bar{\hspace*{-1pt}\bar\cV}_0(\xi)
\ge v_0(\xi).\vadjust{\goodbreak}
\end{equation}

\subsection{The main results}

In the present modified setting, we have the equality between the
second order target problem and its first relaxation. For this, the
following technical condition is needed.

\begin{assum}
\label{assum-HF}
For any $\e>0$, there is an $\dbF$-progressively measurable $\e
$-maximizer $\gamma^{\e}:=\gamma^\e_t(y,z)$ of (\ref{F})
such that, for every $\d>0$,
\begin{eqnarray*}
|\gamma^{\e}_t(y,z)|\le C_{\e,\delta}(1+|y|+|z|),\qquad
\hat\dbP\mbox{-a.s. on } \{\hat a_t \ge\d I_d\},
\mbox{ for some } C_{\e,\d}>0.
\end{eqnarray*}
\end{assum}

Similar to (\ref{Lambda}), for each $i\ge1$, define
%
\begin{eqnarray}
\label{Lambdai}
\L^i &:=& {\esup\limits_{0\le t\le1}} ^{ \dbP^i}\L^i_t,
\nonumber
\\[-8pt]
\\[-8pt]
\nonumber
 \L^i_t
&:= &{\esup\limits_{\dbP\in\cP^\cA_H(t, \dbP^i)}} {}^{\hspace*{-7pt} \dbP^i}
\biggl(\dbE^\dbP_t\biggl[|\xi|^2 + \int_t^1 |\hat F_s (0)|^2
\biggr]\biggr)^{1/2},
\end{eqnarray}
where, as in Proposition~\ref{prop-sup=esssup},
%
\begin{equation}
\label{cPtP0}
\cP^\cA_H(t, \dbP^i) := \{\dbP\in\cP^\cA_H\dvtx  \dbP= \dbP^i
\mbox{ on } \cF_t\}.
\end{equation}
%

\begin{thmm}\label{thm-BankBaum}
Let Assumptions~\ref{assum-H},~\ref{assum-F} and~\ref{assum-HF} hold
true. Assume further that
%
\begin{equation}
\label{Lambdai-integrable}
\dbE^{\dbP^i}[|\L^i|^2] <\infty \qquad\mbox{for all } i\ge1.
\end{equation}
Then for any $\xi\in\hat\dbL^2_0(\dbR)$, we have $\cV_0(\xi
)={\bar\cV}_0(\xi)=\bar{\hspace*{-1pt}\bar\cV}_0(\xi)= v_0(\xi)$, and
existence holds for the problem $\bar{\hspace*{-1pt}\bar\cV}_0(\xi)$. Moreover,
if $F$ has a progressively measurable optimizer, existence also holds
for the problem $\bar\cV_0(\xi)$.
\end{thmm}

This main result $\cV_0(\xi)={\bar\cV}_0(\xi)$ will be proved in
the next subsection. The equality ${\bar\cV}_0(\xi)=\bar{\hspace*{-1pt}\bar\cV
}_0(\xi)$ was already stated in (\ref{VgebarVgehatV0}). The remaining
statements are analogous to the proof of Theorem~\ref{thm-duality}. We
thus omit the proof and only comment on it:
\begin{itemize}
\item We first define for every $i\ge1$ the dynamic problem:
%
\begin{equation}
\label{Vi}
V^i_t := {\esup\limits_{\dbP\in\cP^\cA_H(t, \dbP^i)}} ^{\hspace*{-6pt} \dbP
^i} \cY^\dbP_t(1,\xi),
\end{equation}
where $\cY^\dbP$ is defined by means of the BSDE (\ref{YZa}) and
$\cP
^\cA_H(t,\dbP^i)$ is given in~(\ref{cPtP0}). In light of Proposition
\ref{prop-sup=esssup}, this is the analogue of the process $V$ in
(\ref{V}), except that this is defined $\dbP^i$-a.s. for every $i\ge
1$. However, using the aggregation Theorem 5.1 in~\cite{STZ09a}, we
can aggregate the family $\{V^i, i\ge1\}$ into a universal process
$V$, that is, $V = V^i$, $\dbP^i$-a.s. for all $i\ge1$.
\item Combining the arguments of Lemma 7.2 in~\cite{STZ09a} and
Proposition~\ref{prop-DPP}, we have the dynamic programming principle:
\[
V_{t_1} = {\esup\limits_{\dbP\in\cP^\cA_H(t, \dbP^i)}} ^{\hspace*{-6pt} \dbP
^i} \cY^{\dbP^i}_{t_1}(t_2, V_{t_2}), \qquad\dbP^i\mbox{-a.s. for
all } i\ge1.
\]
\item Exploiting the connection with reflected BSDEs, we then obtain
the decomposition (\ref{decomposition V}) under each $\dbP^i$, and we
conclude by the definition of the problem~$\bar{\hspace*{-1pt}\bar\cV}_0(\xi)$.
\end{itemize}

Our final result shows that, except for the initial second order target
problem, under certain conditions all other problems are not altered by
the modification of this section:

\begin{thmm}\label{thm-regularity}
Let Assumptions~\ref{assum-H},~\ref{assum-F},~\ref{assum-xi} hold,
and assume further that:
\begin{itemize}[--]
\item[--]$F$ is uniformly continuous in $a$ for $a\in D_{F_t}$, and for all
$(t,\omega,y,z)$ and all \mbox{$a\in D_{F_t}$}:
%
\begin{eqnarray}
\label{growth}
|\xi(\omega)|&\le& C(1+\|\omega\|_1)\quad \mbox{and}
\nonumber
\\[-8pt]
\\[-8pt]
\nonumber
|F_t(\omega, y, z, a)|&\le& C(1+\|\omega\|_t+
|y|+|z|+|a^{1\slash2}|),
\end{eqnarray}
\item[--] $\cP^\cA_H$ is dense in $\cP_H$ in the sense that for any $\dbP
=\dbP^{\a} \in\cP_H$ and any $\e>0$:
%
\begin{eqnarray}
\label{dense}
\dbE^{\dbP_0}\biggl[\int_0^1 |(\a^\e_t)^{1\slash2} - \a
_t^{1\slash2}|^2 \,dt\biggr]\le\e
\qquad\mbox{for some }
\dbP^\e= \dbP^{\a^\e}\in\cP^\cA_H.
\end{eqnarray}
Then $v_0(\xi)=v(\xi)$ and, thus, $v_0(\xi)$ is independent from the
choice of the sets~$\cA_0$ and $T_0$.

Assume further that Assumption~\ref{assum-HF} and (\ref
{Lambda-integrability}) hold. Then
\[
{\cV}_0(\xi)=\bar{\cV}_0(\xi)=\bar{\hspace*{-1pt}\bar\cV}_0(\xi)= v_0(\xi
)=v(\xi)=\bar{\hspace*{-1pt}\bar\cV}(\xi) = \bar\cV(\xi).
\]
\end{itemize}
\end{thmm}

\begin{pf} By (\ref{VgebarVgehatV}), Theorems~\ref{thm-duality} and \ref
{thm-BankBaum}, clearly it suffices to prove the first statement.
Since $\cP^\cA_H\subset\cP_H$, we have $v_0(\xi)\le v(\xi)$. Now
for any $\dbP= \dbP^\a\in\cP_H$ and any $\e>0$, let $\dbP^\e=
\dbP^{\a^\e}\in\cP_H^\cA$ satisfy (\ref{dense}). Recall the
$W^\dbP$ defined in (\ref{WP}). Notice that
\begin{eqnarray*}
Y^\dbP_t = \xi(B_\cd) + \int_t^1 F_s(B_\cd,Y^\dbP_s,Z^\dbP
_s,\hat a_s)\,ds - \int_t^1 Z^\dbP_s \hat a_s^{1\slash2}\,dW^\dbP_s,
\qquad0\le t\le1, \dbP\mbox{-a.s.}
\end{eqnarray*}
Let $(\tilde Y^\dbP, \tilde Z^\dbP)$ denote the solution to the
following BSDE under $\dbP_0$:
\begin{eqnarray}
\tilde Y^\dbP_t =
\xi(X^\a_\cd) + \int_t^1 F_s(X^\a_\cd,\tilde Y^\dbP_s,\tilde
Z^\dbP_s, \a_s)\,ds
- \int_t^1 \tilde Z^\dbP_s \a_s^{1\slash2}\,dB_s,\nonumber \\
\eqntext{ 0\le t\le1, \dbP
_0\mbox{-a.s.}}
\end{eqnarray}
By Lemma~\ref{lem-atilde}, the $\dbP$-distribution of $Y^\dbP$ is
equal to the $\dbP_0$-distribution of $\tilde Y^\dbP$. This, together
with the Blumenthal zero--one law, implies that $Y^\dbP_0 = \tilde
Y^\dbP_0$. Similarly, $Y^{\dbP^\e}_0 = \tilde Y^{\dbP^\e}_0$,
where $(\tilde Y^{\dbP^\e}_0,\tilde Z^{\dbP^\e}_0)$ is the solution of
\begin{eqnarray}
\tilde Y^{\dbP^\e}_t =
\xi(X^{\a^\e}_\cd) + \int_t^1 F_s(X^{\a^\e}_\cd,\tilde Y^{\dbP
^\e}_s,\tilde Z^{\dbP^\e}_s, \a^\e_s)\,ds
- \int_t^1 \tilde Z^{\dbP^\e}_s (\a^\e_s)^{1\slash2}\,dB_s,\nonumber\\
\eqntext{ 0\le
t\le1, \dbP_0\mbox{-a.s.}}
\end{eqnarray}
By Proposition 2.1 from El Karoui, Peng and Quenez~\cite{EPQ}, we
deduce that
\begin{eqnarray*}
|Y_0^\dbP-Y_0^{\dbP^\e}|^2 &=& |\tilde Y_0^\dbP-\tilde Y_0^{\dbP
^\e}|^2\\
&\le&
C\dbE^{\dbP_0}\biggl[|\xi(X^{\a}_\cd)-\xi(X^{\a^\e}_\cd
)|^2
+ \int_0^1 |F_t(X^{\a}_\cd,\tilde Y^\dbP_t,\tilde
Z^\dbP_t,\a_t)
\\
&&\hspace*{86pt}\qquad{}-
F_t(X^{\a^\e}_\cd,\tilde Y^{\dbP}_t,\tilde Z^{\dbP}_t,\a^\e_t)
|^2\,dt
\biggr].
\end{eqnarray*}
By (\ref{growth}) we have
\begin{eqnarray*}
|\xi(X^{\a^\e}_\cd)|&\le& C\|X^{\a^\e}\|_1 \le C\|X^\a\|_1 + C\|
X^\a-X^{\a^\e}\|_1;\\
 |F_t(X^{\a^\e}_\cd,\tilde Y^{\dbP}_t,\tilde Z^{\dbP}_t,\a^\e
_t)|&\le& C(\|X^{\a^\e}\|_t+|\tilde Y^{\dbP}_t|+|\tilde Z^{\dbP
}_t|+|\a^\e_t|^{1\slash2})\\
&\le& C(\|X^{\a}\|_1+|\tilde Y^{\dbP}_t|+|\tilde Z^{\dbP
}_t|+|\a_t|^{1\slash2})\\
&&{} + C(\|X^{\a^\e}-X^\a\|_1+|\a^\e
_t-\a_t|^{1\slash2}).
\end{eqnarray*}
It follows from (\ref{dense}) that $\dbE^{\dbP_0}[\sup_{0\le
t\le1}|X^\a_t - X^{\a^\e}_t|^2]\le\e$. Then $|\xi(X^{\a^\e
}_\cd)|^2$ is uniformly integrable under $\dbP_0$ and $|F_t(X^{\a^\e
}_\cd,\tilde Y^{\dbP}_t,\tilde Z^{\dbP}_t,\a^\e_t)|^2$ is
uniformly integrable under $dt\times d\dbP_0$. Now by the uniform
continuity of $\xi$ and $F$ we get $\lim_{\e\to0}|Y_0^\dbP
-Y_0^{\dbP^\e}|=0$. This implies that $Y^\dbP_0 \le v_0(\xi)$ for
all $\dbP\in\cP_H$, and, therefore, $v(\xi)\le v_0(\xi)$.
\end{pf}

A sufficient condition for the uniform continuity of $F$ in terms of
$a$ is that $D_H$ is bounded. We next provide a sufficient condition
for the density condition (\ref{dense}).

\begin{prop}
\label{prop-dense}
Let Assumption~\ref{assum-H} hold and suppose that the domain $D_F$ of
$F$ is independent of $t$. Assume further that $T_0$ is dense in
$[0,1]$, and there exists a countable dense subset $A \subset D_F$ such
that, for all $a\in A$, the constant mapping $a$ is in $\cA_0$. Then
$\cP^\cA_H$ is dense in $\cP_H$ in the sense of (\ref{dense}).
\end{prop}

\begin{pf}(i) We first prove that $\dbP^\a\in\cP^\cA_H$ for any $\a$
taking the following form:
%
\begin{equation}\label{a}
\begin{tabular}{p{300pt}@{}}
There exist $0=t_0<\cds<t_n=1$ in $T_0$ and a finite subset
$A_n \subset A$ s.t.
$\a= \displaystyle\sum_{i=0}^{n-1} \a_{t_i} \1_{[t_i, t_{i+1})} + \a_{t_n}\1_{\{
t_n\}}$ and $\a$ takes values in $A_n$.
\end{tabular}
\end{equation}
In fact, since $A_n\subset\dbS^{>0}_d$ is finite, then $\a$ has both
lower (away from $0$) and upper bounds, and thus $\dbP^\a$ is well
defined. Using the notation in Lemma~\ref{lem-atilde}, we set $a:= \a
\circ\b_\a$. Clearly, $a = \sum_{i=0}^{n-1} a_{t_i} \1_{[t_i,
t_{i+1})}+a_{t_n}\1_{\{t_n\}}$ and $a$ also takes values in $A_n$. By
Lemma~\ref{lem-atilde} we know $\hat a = a$, $dt\times d\dbP^\a
$-a.s. and $\dbP^\a$ satisfies (\ref{ellipticity}). Then it follows
from (\ref{PA}) that $\dbP^\a\in\cP(\cA)$. Moreover, by numerating
$A_n = \{a^i, i=1,\ldots, n\}$, we have $a = \sum_{i=1}^n a^i \1
_{E_i}$, where $E_i := \{\omega\dvtx  a_t(\omega) = a^i, 0\le t \le1\}
$, $i=1,\ldots, n$, form a partition of $\cF_1$. By Lemma 5.2 in \cite
{STZ09a}, we know $\dbP^\a= \dbP^{a^i}$ on $E_i$, that is, $\dbP^\a
(E\cap E_i) = \dbP^{a^i}(E\cap E_i)$ for all $E\in\cF_1$. Since each
$\dbP^{a^i}\in\cP^\cA_H$ satisfies (\ref{cP}), then so does $\dbP
^\a$. This implies that $\dbP^\a\in\cP_H$, and, therefore,
$\dbP^\a\in\cP^\cA_H$.\vspace*{-6pt}
\begin{longlist}[(ii)]
\item[(ii)] Now fix $\dbP^\a\in\cP_H$. Since $\hat a\in D_F$, $dt\times
d\dbP^\a$-a.s. by Lemma~\ref{lem-atilde} we know $\a\in D_F$,
$dt\times d\dbP_0$-a.s. For any $\e>0$, since $\dbE^{\dbP_0}
[\int_0^1 |\a_t|^2 \,dt]<\infty$, by standard arguments there
exists $\dbF$-progressive measurable {{c\`{a}dl\`{a}g}} process $\a
^\e$ such
that $\a^\e$ takes values in $D_F$ and $\dbE^{\dbP_0}[\int
_0^1 |(\a^\e_t)^{1\slash2}-(\a_t)^{1\slash2}|^2 \,dt]\le\e$.
Now by the dense property of $T_0$ and $A$, there exists $\tilde\a^\e
$ in the form (\ref{a}) such that $\dbE^{\dbP_0}[\int_0^1
|(\tilde\a^\e_t)^{1\slash2}-(\a^\e_t)^{1\slash2}|^2 \,dt]\le
\e$. Then $\dbE^{\dbP_0}[\int_0^1 |(\tilde\a^\e_t)^{1\slash
2}-(\a_t)^{1\slash2}|^2 \,dt]\le C\e$. Since $\dbP^{\tilde\a
^\e}\in\cP^\cA_H$ by the above (i), the proof is complete.\quad\qed
\end{longlist}
\noqed\end{pf}

\subsection{\texorpdfstring{Proof of Theorem \protect\ref{thm-BankBaum} [$\cV_0(\xi) = \bar
\cV_0(\xi)$]}{Proof of Theorem 5.2 [V 0(xi) = V 0(xi)]}}

The proof requires the following extension of Bank and Baum \cite
{bankbaum} to the nonlinear case.

\begin{lem}\label{lem-BankBaum}
Let $h_t(\omega,x,z)\dvtx  [0,1]\times\O\times\dbR\times\dbR^d\to
\dbR$ be $\dbF$-progressively measurable, uniformly Lipschitz
continuous in $(x,z)$, and $h(0,0)\in\hat\dbH^2_0(\dbR)$.
For a process $Z\in\hat\dbH^2_0(\dbR^d)$, let $X^Z\in\hat\dbS
^2_0(\dbR)$ denote the aggregating process of the solutions to the
following ODE (with random coefficients) under each~$\dbP^i$:
\[
X^Z_t = x + \int_0^t h_s(X^Z_s, Z_s) \,ds + \int_0^t Z_s \,dB_s, \qquad 0\le
t\le1, \hP\mbox{-a.s.}
\]
Then for any $\e>0$, there exists $Z^\e\in\hat\dbH^2_0(\dbR^d)$
with finite variation, $\hP$-a.s. such that
\[
\sup_{0\le t\le1} |X^{Z^\e}_t - X^Z_t| \le\e,\qquad
\hP\mbox{-a.s.}
\]
\end{lem}

\begin{pf} Recall (\ref{ellipticity}). For $i\ge1$, let $C_i =
C_i(\underline a_{\dbP^i}, \overline a_{\dbP^i})\ge1$ be some
constants which will be specified later. Note that (\ref{ellipticity})
implies $\hat a^{1\slash2} Z\in\dbH^2(\dbP^i, \dbR^d)$. Define
$\tilde\dbP:= \sum_{i=1}^\infty\n_i \dbP^i$, where $\n_1
:=1-\sum_{i=2}^\infty\n_i>0$, and
%
\begin{eqnarray}\label{alphai}
\frac{1}{\n_i}
:=
2^iC_i
\biggl[1+\dbE^i\biggl\{\sup_{0\le t\le1} |X^Z_t|^2 + \int_0^1
[|Z_t|^2 + |\hat a^{1\slash2}_t Z_t|^2 +
|h_t(0,0)|^2]\,dt\biggr\}\biggr]
\nonumber\hspace*{-35pt}
\\[-4pt]
\\[-12pt]
\eqntext{\mbox{for }i\ge2.}
\end{eqnarray}
Then $\tilde\dbP$ is probability measure equivalent to $\hP$, $\dbP
^i \le\n_i^{-1} \tilde\dbP$,
and
%
\begin{equation}
\label{Zintegrability}
X^Z\in\dbS^2(\tilde\dbP, \dbR)\quad\mbox{and}\quad Z, \hat a^{1\slash
2} Z\in\dbH^2(\tilde\dbP, \dbR^d).
\end{equation}
Obviously, it suffices to
find $Z^\e\in\dbH^2(\tilde\dbP, \dbR)$ such that
\[
Z^\e\mbox{ has finite variation and } \sup_{0\le t\le1}|X^{Z^\e
}_t-X^Z_t|\le\e, \qquad\tilde\dbP\mbox{-a.s.}
\]

(1) Denote $X := X^Z$.
As in Bank and Baum~\cite{bankbaum}, we first prove that, for any
$\dbF$-stopping
time $\t$ and any $\tilde X_\t, \tilde Z_\t\in\dbL^2(\tilde\dbP
,\cF_\t)$, there exists
a process $Z^{\e,\t}\in\dbH^2(\tilde\dbP, \dbR^d)$ such that
$Z^{\e,\t}_\t=\tilde Z_\t$, $Z^{\e,\t}$ is absolutely continuous
in $t$ with finite variation on $[\t,1]$, and
%
\begin{equation}\label{Yet}
\tilde\dbP\Bigl[\sup_{\t\le t\le1} e^{-L(t-\t)}|X^{\e,\t}_t - X_t|
\ge
\e+ |\tilde X_\t- X_{\t}|\Bigr]
\le
\e,
\end{equation}
where $L$ is the
uniform Lipschitz constant of $h$ with respect to $x$, and
%
\begin{equation}\label{odetau}
\qquad X_t^{\e,\t} = \tilde X_\t+ \int_\t^t h_s(X^{\e,\t}_s, Z^{\e,\t
}_s) \,ds + \int_\t^t Z^{\e,\t}_s \,dB_s,\qquad
t\ge\t,
\tilde\dbP\mbox{-a.s.}
\end{equation}

For simplicity we assume $\t=0$ and $\tilde X_\t=\tilde x, \tilde
Z_\t=\tilde z$. Set $Z_t:=\tilde z$ for $t<0$, and define $Z^n_t := n
\int_{t-{1}/{n}}^t Z_s \,ds$ for every $n\ge1$. Then $Z^n_0 =
\tilde z$, $Z^n$ is continuous in $t$ with finite variation, and, by
(\ref{Zintegrability}),
\[
\lim_{n\to\infty}\dbE^{\tilde\dbP}\biggl\{\int_0^1 [|Z^n_t -
Z_t|^2 +|\hat a_t^{1\slash2}(Z^n_t - Z_t)|^2]\,dt\biggr\}=0.
\]
Let $X^n$ and $\tilde X$ be defined by $X^n_0=\tilde X_0=\tilde x$ and
\[
dX^n_t = h_t(X^n_t, Z^n_t) \,dt + Z^n_t \,dB_t,
\qquad
d\tilde X_t = h_t(\tilde X_t, Z_t) \,dt + Z_t \,dB_t.
\]
By the Lipschitz property of $h$, it follows from standard estimates on
SDEs that
\begin{eqnarray*}
\lim_{n\to\infty}\dbE^{\tilde\dbP}\Bigl\{\sup_{0\le t\le1}
|X^n_t - \tilde X_t|^2 \,dt\Bigr\}=0\quad
\mbox{and}\quad
e^{-Lt}|\tilde X_t-X_t|\le|\tilde x-x|.
\end{eqnarray*}
Then, for any $\e>0$,
\begin{eqnarray*}
&&\tilde\dbP\Bigl[\sup_{0\le t\le1}e^{-Lt}|X^n_t-X_t| \ge\e+
|\tilde x-x|\Bigr]\\
&&\qquad\le
\tilde\dbP\Bigl[\sup_{0\le t\le1}e^{-Lt}|X^n_t-\tilde X_t| \ge\e
\Bigr]
\\
&&\qquad\le
\tilde\dbP\Bigl[\sup_{0\le t\le1}|X^n_t-\tilde X_t| \ge\e\Bigr]
\longrightarrow0\qquad \mbox{as } n\to\infty.
\end{eqnarray*}
By setting $Z^{\e,\t} := Z^n$ for $n$ large enough so that the above
probability is less than~$\e$, we complete the proof of (\ref{Yet}).
By our construction, notice that
%
\begin{equation}
\label{Zintegrability2}
Z^{\e,\t}_{\t'} \in\dbL^2(\tilde\dbP, \cF_{\t'}) \qquad\mbox{for
every $\dbF$-stopping time } \t'\ge\t.
\end{equation}

(2) In this step, we construct a sequence of $\dbF
$-stopping times $(\t_i)_{i\ge0}$ which yields the required
approximation $(X^\e,Z^\e)$. We initialize our construction by $\t
_0:=0$, $\tilde X_0=X_0$ and $\tilde Z_0$ arbitrary. Let $\e>0$ be
fixed, and set $\e_n := 2^{-n}e^{-L}\e$.

Assume $\t_i$ is defined and $(X^\e, Z^\e)$ have been defined
over $[0,\t_i]$ with $Z^\e_{\t_i}\in\dbL^2(\tilde\dbP, \cF_{\t
_i})$. By (\ref{Yet}) there exists $\tilde Z^{i+1}\in
\dbH^2(\tilde\dbP, \dbR^d)$ which is absolutely continuous in $t$ and
has finite variation on $[\t_i,1]$ such that $\tilde Z^{i+1}_{\t
_i}=Z^\e_{\t_i}$ and
\[
\tilde\dbP\Bigl\{\sup_{\t_i\le t\le1} e^{-L(t-\t_i)}|\tilde
X^{i+1}_t - X_t|
\ge\e_{i+1} + |X^\e_{\t_i}-X_{\t_i}|\Bigr\} \le\e_{i+1},
\]
where $\{\tilde X_t^{i+1},t\in[\t_i,1]\}$ is the solution of the ODE
(\ref{odetau}) with initial condition $\tilde X_{\t_i}^{i+1}=X^\eps
_{\tau_i}$.
Denote
\[
\t_{i+1} := 1\wedge\inf\bigl\{t\ge\t_i\dvtx  e^{-L(t-\t_i)}|\tilde
X^{i+1}_t - X_t| = \e_{i+1} + |X^\e_{\t_i}-X_{\t_i}|\bigr\},
\]
and define
\[
X^\e_t := \tilde X^{i+1}_t,\qquad Z^\e_t:= \tilde Z^{i+1}_t,\qquad \forall
t\in(\t_i,\t_{i+1}].
\]
In particular, it follows from (\ref{Zintegrability2}) that $Z^{\e
}_{\t_{i+1}} \in\dbL^2(\tilde\dbP, \cF_{\t_{i+1}})$.

We remark that, although the filtration $\dbF$ is not right
continuous, since $\tilde X^{i+1}_t - X_t$ is continuous, the $\t
_{i+1}$ defined here is an $\dbF$-stopping time.
Since $\sum_{i=1}^\infty\tilde\dbP(\t_i<1) \le\sum_{i=1}^\infty
\e_i
<1$, it follows from the Borel--Cantelli Lemma that $\tilde\dbP(\t_i<1,
\forall i) = 0$. That is, $(X^\e, Z^\e)$ is well defined on $[0,1]$
and $Z^\e$ is absolutely continuous in $t$ and has finite variation on $[0,1]$.
Moreover, for $t\in[\t_i,\t_{i+1}]$,
\[
\sup_{\t_i\le t\le\t_{i+1}} e^{-L(t-\t_i)}|X^\e_t - X_t|
\le
\e_{i+1} + |X^\e_{\t_i}-X_{\t_i}|.
\]
Then
\begin{eqnarray*}
\sup_{\t_i\le t\le\t_{i+1}} e^{-Lt}|X^\e_t - X_t|
&\le&
e^{-L\t_i}\e_{i+1} + e^{-L\t_i}|X^\e_{\t_i}-X_{\t_i}|\\
&\le&
\e_{i+1} + e^{-L\t_i}|X^\e_{\t_i}-X_{\t_i}|.
\end{eqnarray*}
By induction one can easily see that $ \sup_{0\le t\le1}
e^{-Lt}|X^\e_t - X_t|\le\sum_{i=1}^\infty\e_{i} = e^{-L}\e$, and then
\[
\sup_{0\le t\le1} |X^\e_t - X_t|
\le
\e,\qquad \tilde\dbP\mbox{-a.s.}
\]
(3) It remains to check that $Z^\e\in\dbH^2(\tilde\dbP,
\dbR^d)$.
For any $i, j\ge1$, note that
\begin{eqnarray*}
X^\e_t = X^\e_{\t_{j}} - \int_t^{\t_{j}} h_s(X^\e_s,Z^\e_s) \,ds -
\int_t^{\t_{j}} Z^\e_s \,dB_s,\qquad t\le\t_j, \dbP^i\mbox{-a.s.}
\end{eqnarray*}
%
By the Lipschitz continuity of $h$ and (\ref{ellipticity}), and
following standard arguments, one can easily see that, for some
constant $C_i\ge1$,
\begin{eqnarray*}
&&\dbE^{\dbP^i}\biggl[\int_0^{\t_j}|Z^\e_t|^2\,dt\biggr]\\
&&\qquad\le
C_i \dbE^{\dbP^i}\biggl[|X^\e_{\t_j}|^2 + \int_0^{\t
_j}|h_t(0,0)|^2\,dt\biggr]\\
&&\qquad\le
C_i \dbE^{\dbP^i}\biggl[\sup_{0\le t\le1}|X_t|^2 + \e^2 + \int
_0^1|h_t(0,0)|^2\,dt\biggr\}\qquad
\mbox{for all } j\ge1.
\end{eqnarray*}
Set $C_i$ in (\ref{alphai}) to be the above constant $C_i$. Then by
sending $j\to\infty$, we get
\begin{eqnarray*}
\dbE^{\dbP^i}\biggl[\int_0^1|Z^\e_t|^2\,dt\biggr]
\le
{1\over2^i\n_i}\qquad \mbox{for all } i\ge2.
\end{eqnarray*}
Then
\begin{eqnarray*}
\dbE^{\tilde\dbP}\biggl[\int_0^1|Z^\e_t|^2\,dt\biggr]
&=&
\sum_{i=1}^\infty\n_i \dbE^{\dbP^i}\biggl[\int_0^1|Z^\e
_t|^2\,dt\biggr]\\
&\le&\n_1 \dbE^{\dbP^1}\biggl[\int_0^1|Z^\e_t|^2\,dt\biggr] + \sum
_{i=2}^\infty2^{-i} <
\infty.
\end{eqnarray*}
This completes the proof.
\end{pf}

\begin{pf*}{Proof of Theorem \protect\ref{thm-BankBaum} [$\cV_0(\xi) = \bar\cV
_0(\xi)$]}
In view of (\ref{VgebarVgehatV0}), we only need to
show that $\cV_0(\xi)\le\bar{\hspace*{-1pt}\bar{\cV}}_0(\xi)$ when $\bar{\hspace*{-1pt}\bar
{\cV}}_0(\xi)<\infty$. For any $\e>0$, there exist
$\bar{\hspace*{-1pt}\bar y}< \bar{\hspace*{-1pt}\bar{\cV}}_0(\xi)+\e$ and $\bar{\hspace*{-1pt}\bar Z}\in
\hat\dbH^2_0(\dbR^d)$ such that the corresponding $\bar{\hspace*{-1pt}\bar
Y}_1:=\bar{\hspace*{-1pt}\bar Y}{}^{y,\bar Z}_1\ge\xi, \hP$-a.s. Set $\bar y :=
\bar{\hspace*{-1pt}\bar y}+\e$ and $\bar Z := \bar{\hspace*{-1pt}\bar Z}$. By Assumption \ref
{assum-HF}, we may find $\bar\G\in\hat\dbG^2_0(D_H)\cap\hat\dbH
^2_0(D_H)$ such that the corresponding $\bar Y_1:=\bar Y^{\bar y,\bar
Z,\bar\G}_1
\ge\xi, \hP$-a.s. Denote for $t\in[0,1]$,
\begin{eqnarray*}
Z^0_t:= \int_0^t \bar\G_s \,dB_s,\qquad \zeta_t := \bar Z_t - Z^0_t,\qquad
Y^0_t:= \int_0^t Z^0_s \,dB_s,\qquad X_t:= \bar Y_t - Y^0_t,
\end{eqnarray*}
and
\begin{eqnarray*}
h_t(\omega,x,z) := \tfrac{1}{2} \hat a_t(\omega)\dvtx \bar\G_t(\omega
) - H_t\bigl(\omega, x+Y^0_t(\omega), z+Z^0_t(\omega), \bar\G
_t(\omega)\bigr).
\end{eqnarray*}
One easily checks that $h$ satisfies the conditions of Lemma \ref
{lem-BankBaum}, and $X = X^\zeta$.
Then, there exists $\zeta^\e\in\hat\dbS^2_0(\dbR^d)$ with finite
variation over $[0,1]$ so that
\[
\sup_{0\le t\le1}|X^{\zeta^\e}_t - X_t|\le\e, \qquad\hP\mbox{-a.s.}
\]
Set $Z^\e:= \zeta^\e+ Z^0$, $Y^\e:= X^{\zeta^\e} + Y^0$, and
observe that $d\langle Z^\eps,B\rangle_t=d\langle Z^0,B\rangle
_t=\bar\G_t\,dt$, $\hP$-a.s. Therefore,\vadjust{\goodbreak} $Z^\eps\in\hSM_0(\dbR^d)$.
Setting $y := \bar y$, one can easily check that $Y^\e$ satisfies
(\ref{Yybar}) for given $(y, Z^\e, \bar\G)$. Notice that (\ref
{Yybar}) coincides with (\ref{Yy}) for given $\bar\G$, we have $Y^\e
= Y^{y,Z^\eps}$. Then
\begin{eqnarray*}
Y^{y,Z^\eps}-\bar Y=X^{\zeta^\e} - X
\quad\mbox{and, thus,}\quad
\sup_{0\le t\le1}|Y^{y,Z^\eps}_t - \bar Y_t|\le\e,\qquad \hP\mbox{-a.s.}
\end{eqnarray*}
Let $L$ denote the Lipschitz constant of $H$ with respect to $y$, and
set $y^\e:=y+e^L\e$.
Then
\[
Y^{y^\e,Z^\eps}_t - Y^{y,Z^\eps}_t = e^L\e+ \int_0^t \lambda_s
(Y^{y^\e,Z^\eps}_s - Y^{y,Z^\eps}_s)\,ds,
\]
where $|\lambda_s|\le L$. This leads to $Y^{y^\e,Z^\eps}_1 -
Y^{y,Z^\eps}_1 = e^L\e e^{\int_0^1 \lambda_tdt} \ge\e$, and, thus,
\begin{eqnarray*}
Y^{y^\e,Z^\eps}_1 \ge Y^{y,Z^\eps}_1 + \e\ge\bar Y_1 \ge\xi,\qquad \hP
\mbox{-a.s.}
\end{eqnarray*}
Therefore, $\cV_0(\xi) \le y+e^L\e\le\bar{\hspace*{-1pt}\bar y}+(1+e^L)\e\le
\bar{\hspace*{-1pt}\bar{\cV}}_0(\xi) + (2+e^L)\e$. Since $\e$ is arbitrary,
this provides the required result.
\end{pf*}


\section{Extension}
\label{sect-extension}

In this section we extend our setting in Section~\ref{sect-2target} by
considering $\overline\cP_S$ instead of $\cP_H$ and by removing the
constraints on the domains of $H$ and $F$. In view of the length of
this paper, we shall only formulate the extended problems heuristically
and will not report the details. However, all the results in this paper
can be extended to this new setting.

Let $H_t(\omega,y,z,\gamma)\dvtx  [0,1]\times\O\times\dbR\times
\dbR^d\times\dbR^{d\times d} \to\dbR\cup\{\infty\}$ be a
measurable mapping, and
\[
F_t(\omega,y,z,a)
:=
\sup_{\gamma\in\dbR^{d\times d}} \biggl\{\frac{1}{2}a\dvtx \gamma-
H_t(\omega,y,z,\gamma)\biggr\}, \qquad a\in\dbS^{>0}_d,
\]
be the corresponding conjugate with respect to $\gamma$ which takes
values in $\dbR\cup\{\infty\}$. We assume $D_{H_t}$, the domain of
$H$ in $\gamma$, is independent of $(y,z)$ and contains $0$, $H$ is
uniformly Lipschitz continuous in $(y,z)$ and lower-semicontinuous in
$\gamma$ for all $\gamma\in D_{H_t}$, and $F$ is measurable. Then
the domain $D_{F_t}$ of $F$ in $a$ is also independent of $(y,z)$, and
$F$ is uniformly Lipschitz continuous in $(y,z)$, for all $a\in D_{F_t}$.

Recall the notation $\hat F^0_t:=\hat F_t(0,0)$, and define the
increasing sequence of $\dbF$-stopping times
%
\begin{equation}\label{hattau}
\qquad\hat\tau_n
: =
1\wedge\inf\biggl\{t\ge0\dvtx  \int_0^t \hat F^0_s\,ds\ge n\biggr\},
\qquad n\ge1; \quad\mbox{and}\quad \hat\tau
:= \lim_{n\to\infty} \hat\t_n.
\end{equation}
Notice that
%
\begin{eqnarray}
\label{Ftau}
\int_0^1 \hat F^0_s\,ds& <&\infty\qquad\mbox{on } \bigcup_{n\ge
1}\{\hat\t_n = 1\}\quad \mbox{and}
\nonumber
\\[-8pt]
\\[-8pt]
\nonumber
\int_0^1 \hat F^0_s\,ds &=&\infty
\qquad\mbox{on } \bigcap_{n\ge
1}\{\hat\t_n < 1\}.
\end{eqnarray}
We shall assume further that
\begin{eqnarray*}
\dbE^\dbP\biggl[\int_0^{\hat\t_n} |\hat F^0_s|^2\,ds\biggr]<\infty
\qquad\mbox{for all }
\dbP\in\overline\cP_S
\mbox{ and } n\ge1.
\end{eqnarray*}
For the present extended setting, we introduce the space
$ \hat\dbL^2(\dbR):= \break\bigcap_{\dbP\in\overline{\cP}_S} \dbL
^2(\dbP,\dbR)$, together with
\begin{eqnarray*}
\hat\dbH^2(\dbR^d) &:=& \bigcap_{\dbP\in\overline{\cP}_S}\dbH
^2_{\mathrm{loc}}(\dbP,\dbR^d) := \bigcap_{\dbP\in\overline{\cP
}_S}\bigcap_{n\ge1}\bigl\{Z\dvtx  Z \1_{[0,\hat\tau_n]}\in\dbH^2(\dbP
, \dbR^d)\bigr\},\\
\hat\dbG^2_H(D_H) &:=& \bigcap_{\dbP\in\overline{\cP}_S}\dbG
^2_{\mathrm{loc}}(\dbP,D_H)\\
&:=& \bigcap_{\dbP\in\overline{\cP}_S}\bigcap_{n\ge1}\biggl\{\G
\dvtx  \biggl({1\over2}\hat a \dvtx  \G- H(0,0,\G)\biggr) \1_{[0,\hat\tau
_n]}\in\dbH^2(\dbP, \dbR)\biggr\},
\end{eqnarray*}
and the corresponding spaces for continuous processes (resp.,
semimartingales): $X\in\hat\dbS^2(\dbR) := \bigcap_{\dbP\in
\overline{\cP}_S}\dbS^2_{\mathrm{loc}}(\dbP, \dbR)$ [resp., $\hSM(\dbR
^d):= \bigcap_{\dbP\in\overline{\cP}_S}\SM_{\mathrm{loc}}(\dbP,\break \dbR
^d)$] iff for every $n\ge1$ and $\dbP\in\overline{\cP}_S$,
$X_{.\wedge\hat\tau_n}\in\dbS^2(\dbP, \dbR)$ [resp., $\SM(\dbP,
\dbR^d)$].

Now given $\xi\in\hat\dbL^2(\dbR)$, the second order stochastic
target problem is defined~by
\[
\cV(\xi):= \inf\{y\dvtx  Y^{y,Z}_1 \ge\xi, \overline\cP
_S\mbox{-q.s. for some } Z \in\hSM(\dbR^d)\},
\]
where $Y:= Y^{y,Z}\in\hat\dbS^2(\dbR)$ is defined by the following
ODE (with random
coefficients):
\begin{eqnarray*}
Y_t &=& y - \int_0^t H_s(Y_s,Z_s,\G_s)\,ds +\int_0^t Z_s\circ dB_s,\qquad
t< \hat\tau, \overline\cP_S\mbox{-q.s.}\\
Y_{\hat\tau} &:=& \lim_{n\to\infty} Y_{\hat\tau_n}\qquad
\mbox{on } \bigcup_{n\ge1}\{\hat\t_n =1\},\\
Y_t &:=&\infty\qquad\mbox{for } t\in[\hat\tau,1]
\mbox{ on } \bigcap_{n\ge1}\{\hat\t_n <1\}.
\end{eqnarray*}
Similarly, the extended relaxed problems are as follows:
\begin{eqnarray*}
\bar\cV(\xi) &:=& \inf\{y\dvtx  \exists(\bar Z, \bar\G) \in
\hat\dbH^2(\dbR^d)\times\hat\dbG^2_H(D_H) \mbox{ s.t. } \bar
Y^{\dbP,y,\bar Z, \bar\G}_1 \ge\xi, \\
&&\hspace*{151pt}{}\dbP\mbox{-a.s. for
all } \dbP\in\overline\cP_S\},\\
\bar{\hspace*{-1pt}\bar\cV}(\xi) &:=& \inf\{y\dvtx  \exists\bar{\hspace*{-1pt}\bar Z} \in
\hat\dbH^2(\dbR^d) \mbox{ s.t. } \bar{\hspace*{-1pt}\bar Y}{}^{\dbP,y,\bar{\hspace*{-1pt}\bar
Z}}_1 \ge\xi, \dbP\mbox{-a.s. for all } \dbP\in\overline\cP
_S\},
\end{eqnarray*}
where $\bar Y^\dbP:= \bar Y^{\dbP,y,\bar Z, \bar\G}$ and $\bar{\hspace*{-1pt}\bar Y} {}^\dbP:=\bar{\hspace*{-1pt}\bar Y}{}^{\dbP,y,\bar{\hspace*{-1pt}\bar Z}}$ are defined by
%
\begin{eqnarray}
\label{Yybar2}
\bar Y_t^\dbP
&=&
y + \int_0^t\biggl (\frac{1}{2}\bar\G_s\: \hat a_s-H_s(\bar Y_s^\dbP
, \bar Z_s,\bar\G_s)\biggr)\,ds
\nonumber\\
&&{}+\int_0^t\bar Z_s\,dB_s,\qquad t< \hat\tau, \dbP\mbox{-a.s.}
\nonumber\\
\bar Y_{\hat\tau}^\dbP
&:=&
\lim_{n\to\infty}\bar Y_{\hat\tau_n}^\dbP
\qquad\mbox{on } \bigcup_{n\ge1}\{\hat\t_n =1\},\nonumber\\
\bar Y_t^\dbP&:=&\infty\qquad\mbox{for } t\in[\hat\tau,1] \mbox{ on } \bigcap_{n\ge1}\{\hat\t_n <1\};
\\
\bar{\hspace*{-1pt}\bar Y}_t {}^\dbP&=& y + \int_0^t \hat F_s(\bar{\hspace*{-1pt}\bar
Y}_s {}^\dbP, \bar{\hspace*{-1pt}\bar Z}_s)\,ds +\int_0^t\bar{\hspace*{-1pt}\bar Z}_s\,dB_s,\qquad
t< \hat\tau, \dbP\mbox{-a.s.}
\nonumber\\
\bar{\hspace*{-1pt}\bar Y}_{\hat\tau} {}^\dbP&:=& \lim_{n\to\infty} \bar{\hspace*{-1pt}\bar Y}_{\hat\tau_n} {}^\dbP
\qquad\mbox{on } \bigcup_{n\ge1}\{\hat\t_n =1\},\nonumber\\
\bar{\hspace*{-1pt}\bar Y}_t {}^\dbP&:=&\infty\qquad\mbox{for } t\in[\hat\tau,1]
\mbox{ on } \bigcap_{n\ge1}\{\hat\t_n <1\}.\nonumber
\end{eqnarray}

Finally, we remark that $\dbP[\bigcup_n \{\hat\t_n=1\}]=1$ for all
$\dbP\in\cP_H$. The dual formulation in this extended setting is the
same as the original $v(\xi)$ defined in (\ref{vxi}). That is, for
dual formulation we still use $\cP_H$, instead of $\overline\cP_S$.
Under certain technical conditions, again we can show that $\bar\cV
(\xi) = \bar{\hspace*{-1pt}\bar\cV}(\xi) = v(\xi)$. Moreover, if we extend the
weaker version in Section~\ref{sectweak2target} analogously, similar
results will still hold.


\begin{appendix}

\section*{Appendix}\label{sectappendix}

In this \hyperref[sectappendix]{Appendix} we prove Lemma~\ref{lem-rcpd} and claim (\ref
{PncPS}). We shall use the notation of Lemma~\ref{lem-atilde}.

\begin{pf*}{Proof of Lemma \protect\ref{lem-rcpd}} ($\dbP^{\t,\omega}\in
\overline\cP^{\t(\omega)}_S$). Let $\dbP= \dbP^\a\in\overline
\cP_S$ be given. We emphasize that we shall consider both the strong
formulation $(B, X^\a)$ under $\dbP_0$ and the weak formulation
$(W^\dbP, B)$ under $\dbP$.
We prove the lemma in four steps.

\textit{Step} 1. We first proceed in the strong formulation. Let
$\tilde
\t$ be an arbitrary $\dbF$-stopping time. We claim that
%
\setcounter{equation}{0}
\begin{equation}
\label{rcpd0}
(\dbP_0)^{\tilde\t, \omega} = \dbP^{\tilde\t(\omega)}_0\qquad
\mbox{for } \dbP_0\mbox{-a.e. } \omega\in\O.
\end{equation}
Since $\int_0^1 |\a_s(\omega)| \,ds<\infty$, $\dbP_0$-a.s. clearly
$\int_{\tilde\t(\omega)}^1|\a^{\tilde\t, \omega}_s(\tilde
\omega)|\,ds<\infty$ for $\dbP_0$-a.e. $\omega\in\O$ and $\dbP
^{\tilde\t(\omega)}_0$-a.e. $\tilde\omega\in\O^{\tilde\t
(\omega)}$. Then
%
\begin{equation}
\label{PainPS1}
\dbP^{\a^{\tilde\t, \omega}} \in\overline\cP^{\tilde\t
(\omega)}_S \qquad\mbox{for }\dbP_0\mbox{-a.e. } \omega\in\O.
\end{equation}

We now prove (\ref{rcpd0}), which amounts to say, for $\dbP_0$-a.e.
$\omega$,
%
\begin{eqnarray}
\label{rcpd1}\qquad
\dbE^{\dbP_0^{\t,\omega}}[\xi] = \dbE^{\dbP_0^{\t(\omega
)}}[\xi] \qquad\mbox{for any bounded $\cF^{\t(\omega
)}_T$-measurable r.v. $\xi$} .
\end{eqnarray}
By standard approximating arguments, it suffices to prove (\ref{rcpd1})
by assuming
\begin{eqnarray}
\xi= e^{\lambda_1 \tilde B^{\t,\omega}_{t_1} + \cds+\lambda_n
\tilde B^{\t,\omega}_{t_n}},\nonumber\\
\eqntext{\mbox{where } \tilde B^{\t,\omega
}_{t} := \omega_t \1_{[0, \t(\omega))}(t) + \bigl[\omega_{\t(\omega
)} + B^{\t(\omega)}_t\bigr]\1_{[\t(\omega), T]}}
\end{eqnarray}
for all rational $0< t_1<\cds<t_n\le T$ and $\lambda_1,\ldots
,\lambda_n\in\dbQ^d$. By the countability of rational numbers, we
may allow the exceptional $\dbP_0$-null set to depend on $\xi$.
Moreover, by backward induction, we may assume without loss of
generality that $n=1$ and $t_n=T$. That is, we want to prove, for any
$\lambda\in\dbQ^d$,
%
\begin{equation}
\label{rcpd2}
\dbE^{\dbP_0^{\t,\omega}}\bigl[e^{\lambda B^{\t(\omega)}_T}\bigr] =
 \dbE^{\dbP_0^{\t(\omega)}}\bigl[e^{\lambda B^{\t(\omega)}_T}\bigr]\qquad
 \mbox{for } \dbP_0\mbox{-a.e. } \omega.
\end{equation}

Note that
\[
\dbE^{\dbP_0^{\t(\omega)}}\bigl[e^{\lambda B^{\t(\omega
)}_T}\bigr] = e^{{|\lambda|^2/2}[T-\t(\omega)]}.
\]
Then, by (\ref{Pto}) and the definition of r.c.p.d., (\ref{rcpd2}) is
equivalent to
%
\begin{eqnarray}
\label{rcpd3}
\dbE^{\dbP_0}\bigl[e^{\lambda[B_T - B_\t] }\eta_\t\bigr] = \dbE
^{\dbP_0}\bigl[e^{ {(|\lambda|^2/2)}[T-\t]}\eta_\t\bigr]
\nonumber
\\[-8pt]
\\[-8pt]
\eqntext{\mbox{where } \eta_t:= \f(B_{s_1\wedge t},\ldots, B_{s_m\wedge t}),}
\end{eqnarray}
for any $0<s_1<\cds<s_m\le T$ and any bounded and smooth function $\f$.

To see (\ref{rcpd3}), we first assume $\t$ takes only finitely many
values, and by otherwise merging the partition points, we assume
without loss of generality that $\t$ takes only values $s_1,\ldots,
s_m$. Then, noting that $B_\cd-B_{s_i}$ is a Brownian motion under
$\dbP_0$,
\begin{eqnarray*}
\dbE^{\dbP_0}\bigl[e^{\lambda[B_T-B_\t] }\eta_\t\bigr] &=& \sum
_{i=1}^m \dbE^{\dbP_0}\bigl[e^{\lambda[B_T-B_{s_i}] }\eta_{s_i}\1
_{\{\t= s_i\}} \bigr] \\
&=& \sum_{i=1}^m \dbE^{\dbP_0}\bigl[ \dbE^{\dbP_0}\bigl[e^{\lambda
[B_T-B_{s_i}] }|\cF_{s_i}\bigr]\eta_{s_i} \1_{\{\t= s_i\}} \bigr]\\
&=& \sum_{i=1}^m \dbE^{\dbP_0}\bigl[ e^{ {(|\lambda|^2/2)}(T-s_i) }\eta_{s_i} \1_{\{\t= s_i\}} \bigr] = \dbE^{\dbP_0}
\bigl[e^{ {(|\lambda|^2/2)}[T-\t]}\eta_\t\bigr].
\end{eqnarray*}
In the general case, we may find stopping times $\t_n \downarrow\t$
such that each $\t_n$ takes finitely many values. Then
\[
\dbE^{\dbP_0}\bigl[e^{\lambda[B_T-B_{\t_n}] }\eta_{\t_n}
\bigr] = \dbE^{\dbP_0}\bigl[e^{{(|\lambda|^2/2)}[T-\t_n]}\eta
_{\t_n}\bigr].
\]
Send $n\to\infty$, and note that $\eta$ is continuous in $t$, then
by the Dominated Convergence Theorem we obtain (\ref{rcpd3}), and hence
prove (\ref{rcpd0}).

\textit{Step} 2. We construct the r.c.p.d. for $\dbP$ in weak
formulation. Define
%
\begin{eqnarray}
\label{tildet}
\tilde\t:= \t\circ X^\a\quad \mbox{and}\quad \tilde\a^{\t, \omega} :=
\a^{\tilde\t, \b_\a(\omega)}.
\end{eqnarray}
%
One can easily see that $\tilde\t$ is also an $\dbF$-stopping time.
By the definition of $\dbP^\a$ and the definition of the mapping $\b
_\a$ in Lemma~\ref{lem-atilde}, we have $\t= \tilde\t\circ\b_\a
$, $\dbP^\a$-a.s. Then it follows from (\ref{PainPS1}) that
%
\begin{equation}
\label{PainPS}
\dbP^{\tilde\a^{\t,\omega}} \in\overline{\cP}^{\t(\omega
)}_S \qquad \mbox{for } \dbP^\a\mbox{-a.e. } \omega\in\O.
\end{equation}
%

\textit{Step} 3. We show that $\dbP^{\t,\omega} = \dbP^{\tilde\a
^{\t,\omega}}$ for $\dbP$-a.e. $\omega\in\O$, by assuming the
following claim which will be proved in Step 4 below:
%
\begin{eqnarray}
\label{strong-rcpd-check}
&&\dbE^{\dbP^\a}[\f(B_{t_1\wedge\t},\ldots, B_{t_n\wedge
\t}) \psi(B_{t_1},\ldots, B_{t_n})]
\nonumber
\\[-8pt]
\\[-8pt]
\nonumber
&&\qquad=\dbE^{\dbP
^\a}[\f(B_{t_1\wedge\t},\ldots, B_{t_n\wedge\t})
\psi_\t]
\end{eqnarray}
for any $0<t_1<\cds<t_n\le1$ and bounded and continuous functions $\f
, \psi$, where
\begin{eqnarray}
\psi_\t(\omega)& :=& \dbE^{\dbP^{\tilde\a^{\t, \omega}}}
\bigl[\psi\bigl(\omega(t_1),\ldots,\omega(t_k), \omega
(t)+B^t_{t_{k+1}},\ldots, \omega(t) + B^t_{t_n}\bigr)\bigr]\nonumber\\
\eqntext{\mbox{for } t:= \t(\omega) \in[t_k, t_{k+1}).}
\end{eqnarray}

Indeed, if (\ref{strong-rcpd-check}) is true, then by the arbitrariness
of $\f$ and $(t_1,\ldots,t_n)$, it follows from the definition of
r.c.p.d. that, for $\dbP^\a$-a.e. $\omega\in\O$ and for $t:= \t
(\omega) \in[t_k, t_{k+1})$,
%
\begin{eqnarray}
\label{psito}
\qquad\psi_\t(\omega) = \dbE^{\dbP^{\t, \omega}}\bigl[\psi
\bigl(\omega(t_1),\ldots,\omega(t_k), \omega(t)+B^t_{t_{k+1}},\ldots,
\omega(t) + B^t_{t_n}\bigr)\bigr].
\end{eqnarray}
We remark that the exceptional $\dbP^\a$-null set is supposed to
depend on $\psi$ and $t_1<\cds<t_n$. However, by standard
approximating arguments, one can easily choose a common null set. That
is, there exists a $\dbP^\a$-null set $E_0$ such that, for any
$\omega\notin E_0$, (\ref{psito}) holds for all $(t_1,\ldots, t_n)$
and all bounded continuous functions $\psi$. This clearly implies
that, for $\omega\notin E_0$,
\begin{eqnarray}
\dbE^{\dbP^{\t,\omega}}[\eta] = \dbE^{\dbP^{\tilde\a^{\t
,\omega}}}[\eta]\nonumber\\
 \eqntext{\mbox{for all bounded and } \cF^{\t(\omega
)}_1\mbox{-measurable random variables } \eta.}
\end{eqnarray}
Then $\dbP^{\t,\omega} = \dbP^{\tilde\a^{\t,\omega}}$, for
$\dbP$-a.e. $\omega\in\O$. This, together with (\ref{PainPS}),
proves that $\dbP^{\t,\omega}\in\cP^{\t(\omega)}_S$, for $\dbP
$-a.e. $\omega\in\O$.

\textit{Step} 4. We now prove (\ref{strong-rcpd-check}). For $t:= \t
(\omega) \in[t_k, t_{k+1})$, by definition of $\dbP^{\tilde\a^{\t
, \omega}}$ we have
\begin{eqnarray*}
\psi_\t(\omega) &=& \dbE^{\dbP^{\t(\omega)}_0}\biggl[\psi
\biggl(\omega(t_1),\ldots,\omega(t_k), \omega(t)+\int_t^{t_{k+1}} \bigl(\a
^{\tilde\t, \b_\a(\omega)}_s\bigr)^{1/2} \,dB^{\t(\omega
)}_s,\ldots,\\
&&{}\hspace*{154pt}\omega(t)+\int_t^{t_{n}} \bigl(\a^{\tilde\t, \b_\a
(\omega)}_s\bigr)^{1/2} \,dB^{\t(\omega)}_s\biggr)\biggr].
\end{eqnarray*}
Then, for each $\omega\in\O$, when $t := \tilde\t(\omega) = \t
(X^\a(\omega))\in[t_k, t_{k+1})$,
\begin{eqnarray*}
&&\psi_\t(X^\a(\omega))\\
 &&\qquad= \dbE^{\dbP^{\tilde\t(\omega
)}_0}\biggl[\psi\biggl(X^\a_{t_1}(\omega),\ldots,X^\a_{t_k}(\omega
), X^\a_t(\omega)+\int_t^{t_{k+1}} (\a^{\tilde\t, \omega
}_s)^{1/2} \,dB^{\tilde\t(\omega)}_s,\ldots, \\
&&\hspace*{178pt}\qquad X^\a_t(\omega)+\int_t^{t_{n}} (\a^{\tilde\t, \omega
}_s)^{1/2} \,dB^{\tilde\t(\omega)}_s\biggr)\biggr];
\end{eqnarray*}
note that $(\dbP_0)^{\t,\omega}$-distribution of $(B^{\tilde\t
(\omega)}, \a^{\tilde\t, \omega}_s(B^{\tilde\t(\omega)})$ is
equal to the $(\dbP_0)^\omega_{\t}$-distribution of $(B_\cd-
B_{\tilde\t(\omega)}, \a^{\tilde\t, \omega}(B_\cd-B_{\tilde
\t(\omega)}))$. Recall (\ref{rcpd0}), and note that, for each
$\omega\in\O$,
\begin{eqnarray*}
\a_s(\omega) = \a\bigl(\omega\otimes_{\tilde\t(\omega)} \omega
^{\tilde\t(\omega)}\bigr) = \a^{\tilde\t, \omega}_s\bigl(\omega
^{\tilde\t(\omega)}\bigr).
\end{eqnarray*}
Then
\begin{eqnarray*}
&&\psi_\t(X^\a(\omega)) \\
&&\qquad= \dbE^{(\dbP_0)^\omega_{\tilde\t
}}\biggl[\psi\biggl(X^\a_{t_1}(\omega),\ldots,X^\a_{t_k}(\omega),
X^\a_t(\omega)+\int_t^{t_{k+1}} (\a_s)^{1/2}(B_\cd)
\,dB_s,\ldots,\\
&&\hspace*{179pt}\qquad  X^\a_t(\omega)+\int_t^{t_{n}} (\a_s)^{1/2}(B_\cd)
\,dB_s\biggr)\biggr]\\
&&\qquad= \dbE^{(\dbP_0)^\omega_{\tilde\t}}[\psi(X^\a
_{t_1},\ldots,X^\a_{t_k}, X^\a_{t_{k+1}},\ldots, X^\a_{t_n})
]\\
&&\qquad= \dbE^{\dbP_0}[\psi(X^\a_{t_1},\ldots,X^\a_{t_k}, X^\a
_{t_{k+1}},\ldots, X^\a_{t_n})|\cF_{\tilde\t}
](\omega),\qquad \dbP_0\mbox{-a.e. } \omega\in\O.
\end{eqnarray*}
Then
\begin{eqnarray*}
&&\dbE^{\dbP^\a}[\f(B_{t_1\wedge\t},\ldots, B_{t_n\wedge
\t}) \psi_\t]\\
&&\qquad= \dbE^{\dbP_0}[\f(X^\a_{t_1\wedge\tilde\t},\ldots,
X^\a_{t_n\wedge\tilde\t}) \psi_{\tilde\t}(X^\a)]\\
&&\qquad= \dbE^{\dbP_0}\bigl[\f(X^\a_{t_1\wedge\tilde\t},\ldots,
X^\a_{t_n\wedge\tilde\t}) \dbE^{\dbP_0}[\psi(X^\a
_{t_1},\ldots,X^\a_{t_k}, X^\a_{t_{k+1}},\ldots, X^\a_{t_n})
|\cF_{\tilde\t}]\bigr]\\
&&\qquad= \dbE^{\dbP_0}[\f(X^\a_{t_1\wedge\tilde\t},\ldots,
X^\a_{t_n\wedge\tilde\t}) \psi(X^\a_{t_1},\ldots,X^\a
_{t_k}, X^\a_{t_{k+1}},\ldots, X^\a_{t_n})]\\
&&\qquad= \dbE^{\dbP^\a}[\f(B_{t_1\wedge\t},\ldots,
B_{t_n\wedge\t}) \psi(B_{t_1},\ldots, B_{t_n})].
\end{eqnarray*}
This proves (\ref{strong-rcpd-check}) and hence the lemma.
\end{pf*}

\begin{pf*}{Proof of claim (\protect\ref{PncPS})} Let $\dbP= \dbP^\a$ and
$\dbP^i_t = \dbP^{\a^i}$ for appropriate $\a$ and $\a^i$,
$i=1,\ldots,n$. Define
\begin{eqnarray*}
\bar\a_s := \a_s\1_{[0,t)}(s) + \Biggl[\sum_{i=1}^n \a^i_s \1
_{E^i_t}(X^\a) + \a_s\1_{\hat E^n_t}(X^\a)\Biggr]\1_{[t,1]}(s).
\end{eqnarray*}
Following similar arguments as in the proof of (\ref
{strong-rcpd-check}), one can easily show that, for any $0<t_1<\cds<t_k
= t<t_{k+1}<\cds<t_n$ and any bounded continuous functions~$\f$ and
$\psi$,
\begin{eqnarray*}
&&\dbE^{\dbP^\a}\Biggl[\f(B_{t_1},\ldots, B_{t_k}) \sum
_{i=1}^n \dbE^{\dbP^{\a^i_t}}[\psi(B_{t_1},\ldots,
B_{t_k}, B_t+B^t_{t_{k+1}},\ldots, B_t + B^t_{t_n})]\1
_{E^i_t}\Biggr]\\
&&\qquad=\dbE^{\dbP^{\bar\a}}[\f(B_{t_1},\ldots, B_{t_k})
\psi(B_{t_1},\ldots, B_{t_n})].
\end{eqnarray*}
Then $\dbP^n = \dbP^{\bar\a}$ and one sees immediately that $\dbP
^n\in\overline{\cP}_S$.

Moreover, since each $\dbP^i_t$ satisfies (\ref{cPt}), one can easily
check that $\dbP^{n}$ satisfies all the requirements in Definition
\ref{defn-cP}, and thus $\dbP^{n}\in\cP_H$.
\end{pf*}
\end{appendix}

\section*{Acknowledgment}
We are very grateful to an anonymous referee
for his/her very careful reading of the original manuscript and many useful
suggestions.

%

%


\printaddresses

\end{document}